\RequirePackage{fix-cm}
\documentclass[twocolumn]{svjour3}          
\smartqed  
\usepackage{graphicx}
\usepackage[utf8]{inputenc}
\usepackage{enumerate}
\usepackage[authoryear,round,longnamesfirst]{natbib}
\usepackage{amsmath}
\usepackage{svg}
\usepackage{epsfig}
\usepackage{bm}
\usepackage{amsfonts}
\usepackage{booktabs}
\usepackage{adjustbox}
\usepackage{amsmath}
\usepackage{amsmath,amssymb,algpseudocode,algorithm}

\usepackage{algorithm}
\usepackage{algpseudocode}
\usepackage{relsize}
\usepackage{soul}
\usepackage{diagbox}
\usepackage{mathtools}
\usepackage{mathrsfs}
\usepackage[mathscr]{eucal}
\usepackage[title]{appendix}
\usepackage[hyphens,spaces,obeyspaces]{url}
\usepackage[inline]{enumitem}
\usepackage{tcolorbox}
\usepackage{subfig}

\usepackage{lscape}
\usepackage{listings}
\definecolor{codegreen}{rgb}{0,0.6,0}
\definecolor{codegray}{rgb}{0.5,0.5,0.5}
\definecolor{codepurple}{rgb}{0.58,0,0.82}
\definecolor{backcolour}{rgb}{0.95,0.95,0.92}
 
\lstdefinestyle{mystyle}{
    backgroundcolor=\color{backcolour},   
    commentstyle=\color{codegreen},
    keywordstyle=\color{magenta},
    numberstyle=\tiny\color{codegray},
    stringstyle=\color{codepurple},
    breakatwhitespace=false,         
    breaklines=true,                 
    captionpos=b,                    
    keepspaces=true,                 
    numbers=left,                    
    numbersep=10pt,                  
    showspaces=false,                
    showstringspaces=false,
    showtabs=false,                  
    tabsize=4
}
 
\definecolor{aurometalsaurus}{rgb}{0.43, 0.5, 0.5}
\definecolor{arsenic}{rgb}{0.23, 0.27, 0.29}

\newcommand{\norm}[1]{\left\lVert #1 \right\rVert}

\begin{document} \sloppy

\title{Implicit bulk-surface filtering method for node-based shape optimization and comparison of explicit and implicit filtering techniques}


\author{Reza Najian Asl \and Kai-Uwe Bletzinger}

\authorrunning{R. Najian Asl \and K.-U. Bletzinger} 

\institute{R. Najian Asl \and K.-U. Bletzinger \at
              Lehrstuhl für Statik, Technische Universität München, Arcisstr. 21, 80333 München, Germany \\
              Tel.: +49-89-28922422\\
              Fax:  +49-89-28922421\\
              \email{reza.najian-asl@tum.de}           
}

\date{Received: date / Accepted: date}

\maketitle

\begin{abstract}
This work studies shape filtering techniques, namely the convolution-based (explicit) and the PDE-based (implicit), and introduces an implicit bulk-surface filtering method to control the boundary smoothness and preserve the internal mesh quality simultaneously in the course of bulk (solid) shape optimization. To that end, volumetric mesh is filtered by the solution of pseudo-solid governing equations which are stiffened by the mesh-Jacobian and endowed with the Robin boundary condition which involves the Laplace-Beltrami operator on the mesh boundaries. Its superior performance from the non-simultaneous (sequential) treatment of boundary and internal meshes is demonstrated for the shape optimization of a complex solid structure. Well-established explicit filters, namely Gaussian and linear, and the Helmholtz/Sobolev-based (implicit) filter are critically examined in terms of consistency (rigid-body-movement production), geometric characteristics and the computational cost. It is shown that the implicit filtering is numerically more efficient and unconditionally consistent, compared to the explicit one. Supported by numerical experiments, a regularized Green's function is introduced as an equivalent explicit form of the Helmholtz/Sobolev filter. Furthermore, we give special attention to derive mesh-independent filtered sensitivities for node-based shape optimization with non-uniform meshes. It is shown that the mesh independent filtering can be achieved by scaling discrete sensitivities with the inverse of the mesh mass matrix.


\end{abstract}

\section{Introduction}
Filtering is now well-established and has proved effective in discrete topology and shape optimization. For a given filter length scale, it should regularize the optimization problem independent of the discretization of the underlying geometry. Together with adjoint-based sensitivity analysis to determine the discrete gradients, the filtering or regularization has become a very successful procedure for large-scale optimization problems in industry. Techniques for filtering are divided in two main categories whose similarities and differences are going to be discussed in this paper. In the first approach, smoothing is applied implicitly by solving elliptic PDEs whose inverse operator is a local smoother, the so-called implicit filtering. In the second category,  raw field is smoothed by its convolution with a kernel function, the so-called explicit filtering. Both techniques are conventionally used to smooth either design variables or design gradients. The former is more consistent with the optimization problem formulation, whereas the later may confuse the optimizer and disturb the convergence, due to discrepancies between the filtered and true (raw) gradients. 

Explicit filtering was introduced by \cite{sigmund1994design} in order to eliminate two well known problems in structural topology optimization, namely the checkerboard problem and the mesh-dependency problem. Since then, the method has become a standard tool in structural design and optimization, mainly due to its simple formulation and robustness. In shape optimization, the explicit filtering has been applied to structural \citep{bletzinger2005computational,bletzinger2010optimal,le2011gradient,firl2012shape} as well as fluid and aerodynamic \citep{stuck2011adjoint,hojjat2014vertex,kroger2015cad,najian2017consistent,ghantasala2021realization,antonau2022latest} applications. Furthermore, \cite{bletzinger2014consistent} has established perfect analogy between the explicit filtering and parameterization techniques used for the shape optimization, especially CAD-based techniques.  

The traction method \citep{azegami1996domain} and the Sobolev-gradient smoothing \citep{Jameson2000} were pioneering works on the PDE-based or implicit filtering. The main idea is to project the gradients into a smoother design space, the Hilbert or Sobolev space, by letting the gradients to be the solution of an elliptic equation. \cite{schmidt2008shape,eppler2009preconditioning} have demonstrated mathematically that the Sobolev smmothing can be interpreted as a reduced shape Hessian operator which turns the steepest descent into a quasi-Newton step. Following this interpretation, \cite{DICK2022105568} have used the Sobolev smoothing as a preconditioner for discrete shape sensitivities to accelerate the convergence in CAD-based shape optimization, even though the parameterization is smooth itself. \cite{mohammadi2009applied} also observed that the conditioning of the optimization problem is more favorable with the Sobolev-based filtering than with the CAD-based filtering. Recently, the implicit filtering has come to the attention of topology optimization community, after the early works by \cite{lazarov2011filters,kawamoto2011heaviside}. The main motivation behind is that, it has better performance profiles compared to the classical or explicit filtering.

This work studies the explicit and implicit filters in the context of node-based shape optimization. They are discussed qualitatively and quantitatively in terms of the consistency between shape gradients and shape updates, interpolation of the domain edges and the sensitivity of filtered fields to the discretization of the underlying geometry, i.e. mesh-dependency.  Furthermore, we present a bulk-surface filtering method which is designed to control the boundary smoothness and prevent the internal mesh distortion simultaneously. This is similar to the enhanced traction method \citep{azegami2006smoothing} which computes shape updates by solving a pseudo-elastic linear problem which is stiffened by springs on the domain boundary and is loaded by discrete shape gradients. In this work, with the aim to maintain mesh quality and reduce the frequency of remeshing during the shape optimization, the so-called Jacobian-based stiffening 
\citep{tezduyar,stein,tonon2021linear} is applied to the pseudo-elastic solid model. Moreover, to ensure smooth boundary geometry shapes, the model is endowed with the Robin boundary condition which involves the Laplace-Beltrami operator on the boundary mesh. Furthermore, we apply the filtering to discrete shape rather than discrete shape sensitivities, to avoid possible discrepancies between the filtered and actual sensitivities. By doing so, a smoother design space, Hilbert or Sobolev space, is introduced in parallel to the geometry space where shape updates are applied. Mathematically speaking, the actual optimization problem is defined in the smooth control space and control sensitivities as well as shape update are derived consistently. Last but not least, since discrete sensitivities are consistent nodal values which show size effects in the gradients, inverse of the mesh mass matrix is used as a preconditioner for the discrete control sensitivities to avoid the mesh dependency issue.

This paper is structured as follows: In Section 2, the explicit shape filtering is revisited and its discrete form is derived. Therein, consistency and mesh-dependency issues are elaborated and the mesh mass matrix preconditioner is introduced. In Section 3, the implicit shape filtering is presented for surface geometries as well as bulk domains using the previously described bulk-surface formulation. Finally, Section 4 compares critically the explicit and implicit filtering techniques regarding  geometric characteristics and the computational cost. Therein, performance of the developed bulk-surface filter is studied for the shape optimization of a complex solid structure and comparisons made against the non-simultaneous (sequential) treatment of boundary and internal meshes.


\section{Explicit shape filtering}
\label{s:explicit-shape-filtering}
This technique generates the three-dimensional geometry at point $ \boldsymbol{ x}_{0} = (x^1_{0},x^2_{0},x^3_{0} ) $ of the design surface $\Gamma$ from the control field $ \boldsymbol{s}=(s^1,s^2,s^3) $ via a convolutional filtering operation:
\begin{equation}
\label{e:conv_filter}
	\boldsymbol{ x}_{0} =  \int\limits_{\Gamma} F(\boldsymbol{ x},\boldsymbol{ x}_{0}) \: \boldsymbol{s}(\boldsymbol{ x}) \: d\Gamma = \int\limits_{\Sigma} F(\text{\bfseries x},\text{\bfseries x}_{0}) \: s \: d\Sigma
\end{equation}

where $F$ is the filter (kernel) function which controls the properties of filtering and subsequently the produced shapes, $\Sigma$ is the portion of  $\Gamma $ which lies within the filter's support span. In free-form shape optimization, unlike topology and CAD-based shape optimization, the control field is unknown \textit{apriori} and it can be calculated by inverse filtering (deconvolution). Gaussian and linear hat functions are the most commonly used kernels in shape and topology optimization and they are defined respectively as 

\begin{subequations}
	\label{e:gauss_hat_kernels}
	\begin{align}
	F^G(\boldsymbol{ x},\boldsymbol{ x}_{0}) & = \frac{1}{r^E_\mathsmaller{\Gamma} \sqrt{2\pi} } \,	\mathlarger{\mathlarger{{e}^{\frac{-1}{2}(\norm{\boldsymbol{ x}-\boldsymbol{ x}_{0}}/r^{E}_\mathsmaller{\Gamma})^2}}} \label{e:gaussian}\\	
	F^L(\boldsymbol{ x},\boldsymbol{ x}_{0}) & = max \left(0, \dfrac{r^{E}_\mathsmaller{\Gamma}-\norm{\boldsymbol{ x}-\boldsymbol{ x}_{0}}}{r^{E}_\mathsmaller{\Gamma}}\right) \label{e:hat}
	\end{align}
\end{subequations}

\begin{figure}
	\includegraphics[keepaspectratio,width=0.5\textwidth]{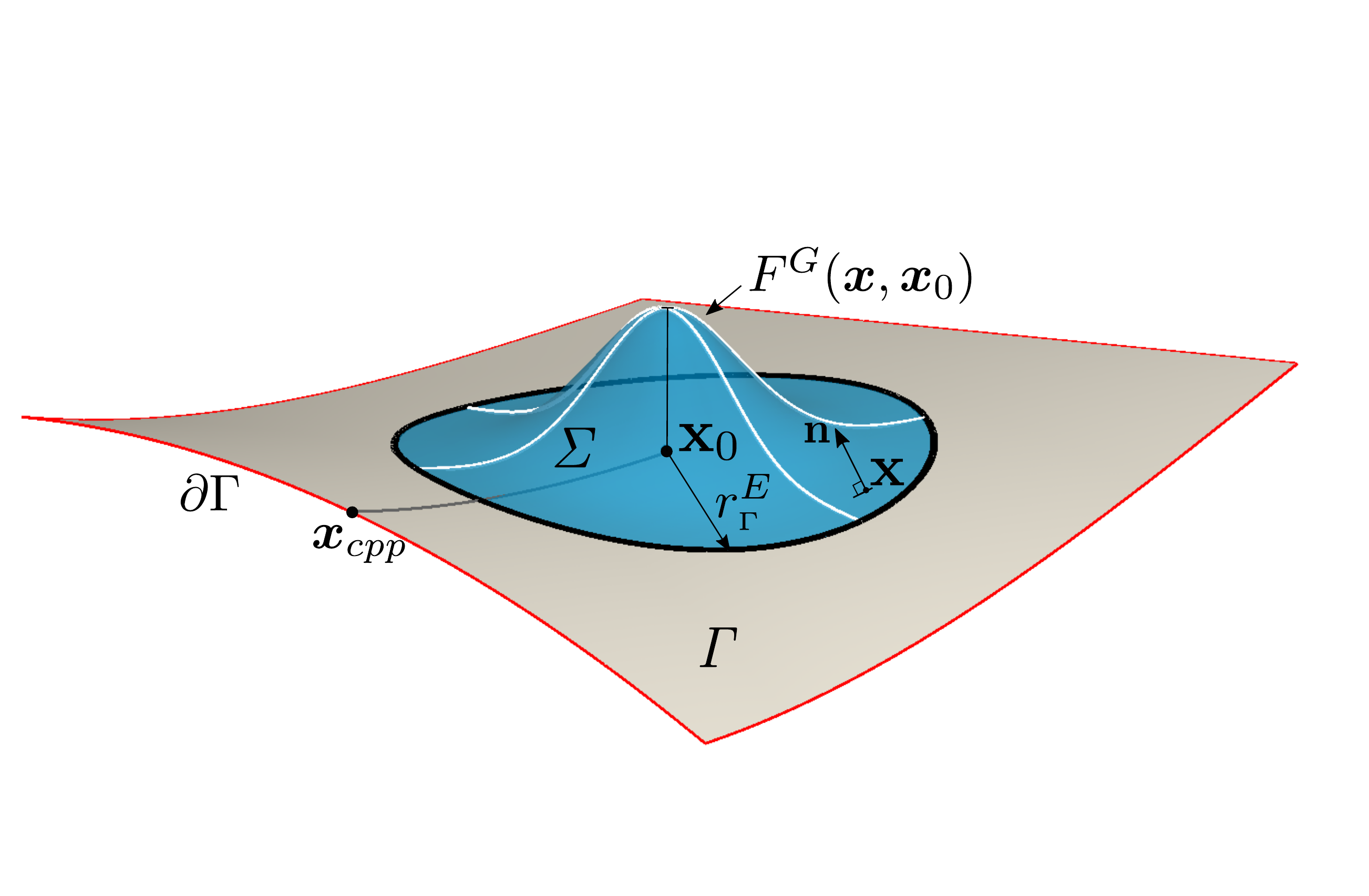}
	\caption{Notional schematic of the design surface ($\Gamma$), the Gaussian filter function ($F$), the integration area ($\Sigma$) and the closest point projection of $\boldsymbol{x}_{0}$ ($\boldsymbol{x}_{cpp}$)} onto the design surface boundary ($\partial \Gamma$).
	\label{fig:VM}      
\end{figure}

$r^{E}_\mathsmaller{\Gamma}$ is the surface filter radius of the explicit filtering; $\boldsymbol{ x}$ is the continuous 3D vector with one component for each spatial coordinate, i.e., $\boldsymbol{ x} = (x^1,x^2,x^3)$; $\norm{\boldsymbol{ x}-\boldsymbol{ x}_{0}}$ is the Euclidean distance to the center of the filter $\boldsymbol{ x}_{0}$. See Fig. \ref{fig:VM} for a notional schematic of the used notation.

In shape optimization, deign surface is often open and limited by adjacent non-design areas, which must not be modified. In other words, the edges of the design surface are fixed or must be interpolated during the optimization. As a matter of fact, edge interpolation is a Dirichlet boundary condition to the optimization problem and must be satisfied strongly and point-wise. In the explicit shape filtering, damping \citep{kroger2015cad,baumgartner2020grid} is used to suppress modification of the fixed boundaries and smoothly transit between design and not-design areas. The implementation of damping is straightforward and one only needs to multiply the filter operator with a damping function which is computed based on the distance to the non-design domain. Then Eq.\ref{e:conv_filter} is updated as 

\begin{subequations}
	\label{e:contA_updated}
	\begin{align}
	& \boldsymbol{x}_{0} = D(\boldsymbol{x}_{0},\boldsymbol{x}_{cpp})\,\int\limits_{\Sigma} F(\boldsymbol{ x},\boldsymbol{ x}_{0}) \: \boldsymbol{s} \: d\Sigma  \\
	& D(\boldsymbol{x}_{0},\boldsymbol{x}_{cpp}) = 1-F(\boldsymbol{ x}_{0},\boldsymbol{ x}_{cpp})  \label{e:contA_sub}
	\end{align}
\end{subequations}

where $D$ is the damping function, $\boldsymbol{x}_{cpp}$ is determined based on the so-called closest point projection (CPP) of $\boldsymbol{x}_{0}$ onto the design surface boundary $\partial \Gamma$.

For numerical analysis, the shape governing equation of the design surface should be discretized. Although one can use different resolutions to discretize the left- and right-hand side of Eq. \ref{e:conv_filter}, in this work the same grid is used for the sake of simplicity. Applications of non-matching grids can be found in \cite{asl2019shape}. We use piecewise linear finite element functions to approximate the geometry as well as data fields, i.e. isoparametric finite elements. Therefore, actual geometry and its control field  within each element are approximated as $	x^i \approx \boldsymbol{N}^i \boldsymbol{x}^i,\;\; s^i \approx \boldsymbol{N}^i \boldsymbol{s}^i; \;\; i \in \{1,2,3\}$, where $\boldsymbol{N}^i$ is the vector of element shape functions in the $i$th Cartesian direction, $\boldsymbol{x}^i$ and $\boldsymbol{s}^i$ are the $i$th component of nodal coordinate vector and nodal control point vector, respectively. Then, the geometry at node $i$ reads

\begin{equation}
\boldsymbol{x}_{\mathsmaller{\Gamma},i} \,=\, D(\boldsymbol{x}_{\mathsmaller{\Gamma},i}) \sum\limits_{k\in\mathbb{E}_{e,i}}\int_{\Gamma^e} F(\boldsymbol{x}_{\mathsmaller{\Gamma},i},\boldsymbol{x}_{\mathsmaller{\Gamma},k})\,\boldsymbol{N}^T_\mathsmaller{\Gamma}\boldsymbol{N}_\mathsmaller{\Gamma} \,\boldsymbol{s}_{\mathsmaller{\Gamma},k}\,d\Gamma
\end{equation}

where $\mathbb{E}_{e,i}$ is the set of elements lying in the filter support domain of node $i$, subscript $\Gamma$ indicates that the quantity belongs to the surface. Note that the filter function is evaluated at the mesh points (not element Gauss points), i.e. it is discretized. In the vector-matrix format, the explicit shape filtering reads
\begin{equation}
\label{e:explicit-vector-matrix}
\begin{aligned}
\boldsymbol{x}_{\mathsmaller{\Gamma}} & = \boldsymbol{D}_\mathsmaller{\Gamma}\,.\, \boldsymbol{W}_\mathsmaller{\Gamma}\,.\,\boldsymbol{M}_\mathsmaller{\Gamma}\,.\, \boldsymbol{s}_\mathsmaller{\Gamma} =  \boldsymbol{A}^E_\mathsmaller{\Gamma} \,.\, \boldsymbol{s}_\mathsmaller{\Gamma} \\
 \boldsymbol{M}_\mathsmaller{\Gamma} &= \sum\int_{\Gamma^e} \boldsymbol{N}^T_\mathsmaller{\Gamma}\boldsymbol{N}_\mathsmaller{\Gamma} d\Gamma 
\end{aligned}
\end{equation}

$\boldsymbol{x}_{\mathsmaller{\Gamma}}\in \mathbb{R} ^{3n_\mathsmaller{\Gamma}\times1}$ is the nodal coordinates vector of the design boundary mesh, $\boldsymbol{A}^E_\mathsmaller{\Gamma}\in \mathbb{R} ^{3n_\mathsmaller{\Gamma}\times3n_\mathsmaller{\Gamma}}$ is the explicit filtering matrix composed of the surface mesh mass matrix $\boldsymbol{M}_\mathsmaller{\Gamma}$, the weighting matrix $\boldsymbol{W}_\mathsmaller{\Gamma}$ whose entries are  $\boldsymbol{W}_\mathsmaller{\Gamma}(i,j)=F(\boldsymbol{x}_{\mathsmaller{\Gamma},i},\boldsymbol{x}_{\mathsmaller{\Gamma},j}) \,.\,\boldsymbol{I}_{3\times3}$ and the diagonal damping matrix $\boldsymbol{D}_\mathsmaller{\Gamma}(i,i)= (1-F(\boldsymbol{x}_{\mathsmaller{\Gamma},i},\boldsymbol{ x}_{cpp})) \,.\,\boldsymbol{I}_{3\times3}$. A drawback of such convolutional filtering is that the \textit{smoothing} effect becomes mesh dependent, especially on unstructured meshes and finite domains with open boundaries. It is clear that the discrete surface $\boldsymbol{x}_{\mathsmaller{\Gamma}}$ is constructed based on a weighted sum of the control values at the mesh points, where the weighting is purely a function of the distances to mesh points. The standard remedy found in literature is to rescale the filter function so that it has a unit volume everywhere in the domain, i.e.     
\begin{equation}
\int\limits_{\Gamma} F(\boldsymbol{ x},\boldsymbol{ x}_{0}) \: d\Gamma =1.0, \;\; \;\;\forall \boldsymbol{x}_{0} \;\; \text{on} \;\; \Gamma   
\end{equation}
This can be achieved easily by dividing the filter function by the inverse of its volume, i.e. $F(\boldsymbol{ x},\boldsymbol{ x}_{0})/=\int\limits_{\Gamma} F(\boldsymbol{ x},\boldsymbol{ x}_{0}) \: d\Gamma$. This scaling has two effects. On one hand, it results in the so-called \textit{consistency} property for the filtering, meaning that a constant variation in the control field $\delta\boldsymbol{s}=[1,1,1]$ results in a uniform update in the geometry  $\delta\boldsymbol{x}=[1,1,1]$. On the other hand, it destroys the self-adjoint property of the filter function and accordingly the symmetry property of the weighting matrix. The consistent explicit filtering matrix reads   
\begin{equation}
\label{e:normalized-explicit-vector-matrix}
\boldsymbol{A}^E_\mathsmaller{\Gamma} = \boldsymbol{D}_\mathsmaller{\Gamma}\,.\, \boldsymbol{V}^{-1}_\mathsmaller{\Gamma}\,.\,\boldsymbol{W}_\mathsmaller{\Gamma}\,.\,\boldsymbol{M}_\mathsmaller{\Gamma}  
\end{equation}
where $\boldsymbol{V}_\mathsmaller{\Gamma}\in \mathbb{R} ^{3n_\mathsmaller{\Gamma}\times3n_\mathsmaller{\Gamma}}$ is a diagonal matrix and its entries are basically the volume under nodal filter functions, i.e. $\boldsymbol{V}_\mathsmaller{\Gamma}(i,i)=  \sum\limits_{k\in\mathbb{E}_{e,i}}\int_{\Gamma^e} F(\boldsymbol{x}_{\mathsmaller{\Gamma},i})\,N_\mathsmaller{\Gamma} \,\,d\Gamma\,.\,\boldsymbol{I}_{3\times3}$. Replacing $\boldsymbol{D}_\mathsmaller{\Gamma}$ by the identity matrix, it is noted that the \textit{consistency} property is reflected as unit row sums in the filtering matrix. Also note that the matrix is a large and dense matrix, however, it can be applied in a matrix-free mode by computing the weights in advance, storing them inside each node and cycling through the nodes. At first glance, lack of the self-adjoint property does not seem relevant, however a closer look at the response sensitivity analysis and shape update rule that will follow, reveals some inconsistency issues.
 
\subsection{Sensitivity analysis and update rule}
In this work, filtering is applied to the discrete shape rather than the discrete shape sensitivities, meaning that parallel to the geometry space $x$ where shape sensitivities are evaluated, the design space $s$ lives where the optimization problem is defined. Applying chain rule of differentiation, control gradients of the response $J$ read 

\begin{equation}
\label{e:sensitivty-mapping}
\frac{dJ}{d\boldsymbol{s}_\mathsmaller{\Gamma}} \, =\, (\boldsymbol{A}^E_\mathsmaller{\Gamma})^T \,.\, \frac{dJ}{d\boldsymbol{x}_\mathsmaller{\Gamma}} = \boldsymbol{M}_\mathsmaller{\Gamma} \,.\, \boldsymbol{W}_\mathsmaller{\Gamma} \,.\, \boldsymbol{V}^{-1}_\mathsmaller{\Gamma} \,.\, \boldsymbol{D}_\mathsmaller{\Gamma}\,.\, \frac{dJ}{d\boldsymbol{x}_\mathsmaller{\Gamma}}
\end{equation}

where $dJ/d\boldsymbol{s}_\mathsmaller{\Gamma}$ and $dJ/d\boldsymbol{x}_\mathsmaller{\Gamma}$ are the discrete control and shape sensitivities of the response respectively. In fluid and structural optimization, response functions are typically a smooth differentiable integral functional whose continuous variation reads
\begin{equation}
\delta J(\boldsymbol{x}) \,=\, \int\limits_{\Gamma} \delta j(\boldsymbol{x}) \, d\Gamma 
\end{equation}
where the spatial variation of the integral is neglected. Then, consistent or discrete gradients can be derived by basis functions of finite elements \citep{oden1972mixed} as
\begin{equation}
\label{e:discrete-continuous-sens}
\frac{dJ}{d\boldsymbol{x}}_\mathsmaller{\Gamma} \,=\, \boldsymbol{M}_\mathsmaller{\Gamma} \,.\, \frac{dj}{d\boldsymbol{x}}_\mathsmaller{\Gamma}
\end{equation} 
$dj/d\boldsymbol{x}_\mathsmaller{\Gamma}$ is the continuous gradient field sampled at mesh points and associated to the discrete/consistent field $dJ/d\boldsymbol{x}_\mathsmaller{\Gamma}$ via the mesh mass matrix $\boldsymbol{M}_\mathsmaller{\Gamma}$. Since the control field $\boldsymbol{s}_\mathsmaller{\Gamma}$ is also a continuous field sampled at mesh points, it must be updated with a field of the same type in the direction of steepest descent. Therefore, inspired by Eq. \ref{e:discrete-continuous-sens} it seems reasonable to develop the following control update rule:  
\begin{equation}
\label{e:quasi-newton-update-rule}
\begin{split}
	\Delta \boldsymbol{s}_\mathsmaller{\Gamma} \, & =\, -\alpha  \,.\, \frac{dj}{d\boldsymbol{s}_\mathsmaller{\Gamma}} \,=\, -\alpha  \,.\, \boldsymbol{M}_{\mathsmaller{\Gamma}}^{-1}  \,.\, \frac{dJ}{d\boldsymbol{s}_\mathsmaller{\Gamma}} \\
		& =\, -\alpha  \,.\, \boldsymbol{M}_{\mathsmaller{\Gamma}}^{-1}  \,.\, \boldsymbol{M}_\mathsmaller{\Gamma} \,.\, \boldsymbol{W}_\mathsmaller{\Gamma} \,.\, \boldsymbol{V}^{-1}_\mathsmaller{\Gamma}  \,.\, \boldsymbol{D}_\mathsmaller{\Gamma}\,.\, \frac{dJ}{d\boldsymbol{x}_\mathsmaller{\Gamma}} \\
		& =\, -\alpha  \,.\, \boldsymbol{W}_\mathsmaller{\Gamma} \,.\, \boldsymbol{V}^{-1}_\mathsmaller{\Gamma} \,.\, \boldsymbol{D}_\mathsmaller{\Gamma} \,.\,  \frac{dJ}{d\boldsymbol{x}}_\mathsmaller{\Gamma} \\
		& =\, -\alpha  \,.\, \boldsymbol{W}_\mathsmaller{\Gamma} \,.\, \boldsymbol{V}^{-1}_\mathsmaller{\Gamma}  \,.\, \boldsymbol{D}_\mathsmaller{\Gamma}\,.\, \boldsymbol{M}_\mathsmaller{\Gamma} \,.\, \frac{dj}{d\boldsymbol{x}}_\mathsmaller{\Gamma}
\end{split}
\end{equation}
where $\alpha$ is a constant step size, $dj/d\boldsymbol{s}_\mathsmaller{\Gamma}$ is the scaled/continuous control sensitivities which are linked to the consistent ones $dJ/d\boldsymbol{s}_\mathsmaller{\Gamma}$ via the inverse mass matrix. One can interpret the proposed update rule  as a quasi-Newton step with diagonal approximation of the Hessian matrix by the scaled mass matrix. A careful look at Eqs. \ref{e:sensitivty-mapping},\ref{e:quasi-newton-update-rule} shows that the desirable consistency property (unit row sum) in the explicit shape filtering (Eqs. \ref{e:explicit-vector-matrix},\ref{e:normalized-explicit-vector-matrix}) does not hold for the sensitivity filtering. This means that a uniform continuous gradient field, e.g. $dj/d\boldsymbol{x}_\mathsmaller{\Gamma}=[1]_{3n_\mathsmaller{\Gamma}\times1}$, does not result in a uniform gradient field in the control space and subsequently a uniform geometry update, i.e. 

\begin{figure}
	\centering
	\includegraphics[keepaspectratio,width=0.4\textwidth]{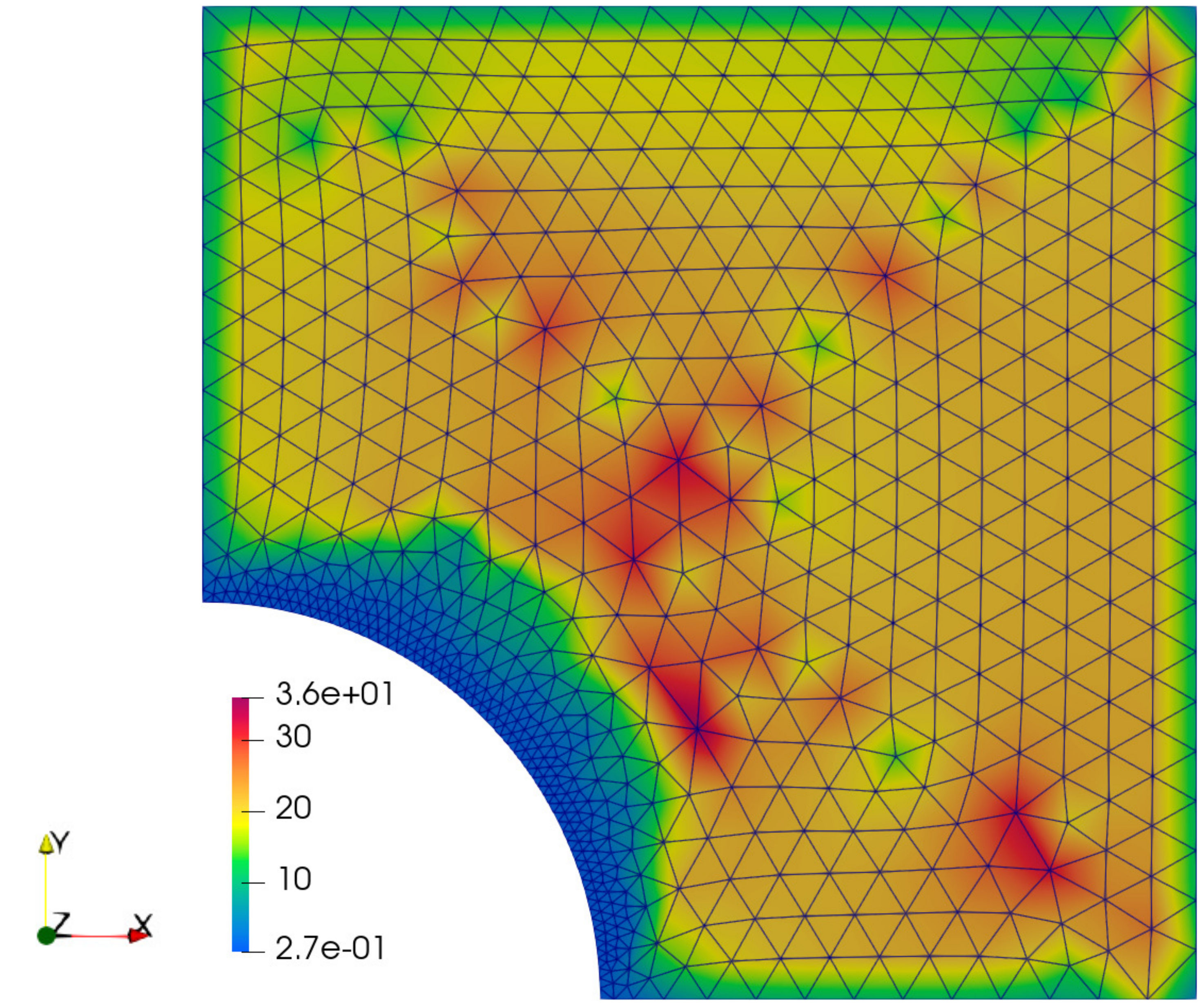}
	\caption{A perforated plate meshed non-uniformly. Consistent/discrete nodal sensitivities $dJ/d\boldsymbol{x}_\mathsmaller{\Gamma}$ associated with the uniform continuous sensitivity field $dj/d\boldsymbol{x}_\mathsmaller{\Gamma}=[(0,0,1)]_{n_\mathsmaller{\Gamma}\times1}$.}
	\label{fig:non-uni-plate}
\end{figure}

\begin{equation}
\label{e:explicit-update-rule}
\begin{split}
\Delta \boldsymbol{x}_\mathsmaller{\Gamma} \, & =\, \boldsymbol{A}^E_\mathsmaller{\Gamma} \,.\, \Delta \boldsymbol{s}_\mathsmaller{\Gamma} \\
 \, & =\,  -\alpha  \,.\, \boldsymbol{A}^E_\mathsmaller{\Gamma} \,.\, \boldsymbol{W}_\mathsmaller{\Gamma} \,.\, \boldsymbol{V}^{-1}_\mathsmaller{\Gamma}  \,.\, \boldsymbol{M}_\mathsmaller{\Gamma} \,.\, [1]_{3n_\mathsmaller{\Gamma}\times1} \\
 \, & \neq\, -\alpha  \,.\, [1]_{3n_\mathsmaller{\Gamma}\times1}
\end{split}
\end{equation}       
where for the sake of simplicity no damping or design boundary interpolation is considered, i.e. $\boldsymbol{D}_\mathsmaller{\Gamma}(i,i)= \boldsymbol{I}_{3\times3}$. The above holds in particular for open design surfaces, independent of the mesh type.
\subsection{Consistency check}
\label{s:explicit-consistency-check}
In the context of gradient-based shape optimization, the \textit{consistency} property means rigid-body-movement production for a given translational (uniform) sensitivity field, independent of the mesh discretization. \cite{kroger2015cad} have discussed consistency aspects of the explicit filtering in great detail. They strictly require the filter to compute constant and linear filtered fields, respectively, for constant and linear distributions of the continuous sensitivity field. In shape and topology optimization, the filter matrix can also be only applied to the design sensitivities rather than the design variables themselves, the so-called sensitivity filtering technique. However, in shape and density filtering, the matrix is used twice; once in the design update rule (Eq. \ref{e:explicit-update-rule}) and once in the sensitivity analysis using a transpose operation (Eq. \ref{e:sensitivty-mapping}). 

\begin{figure}
	\centering
	\begin{tabular}{@{}c@{}}
		\includegraphics[width=0.7\linewidth,keepaspectratio]{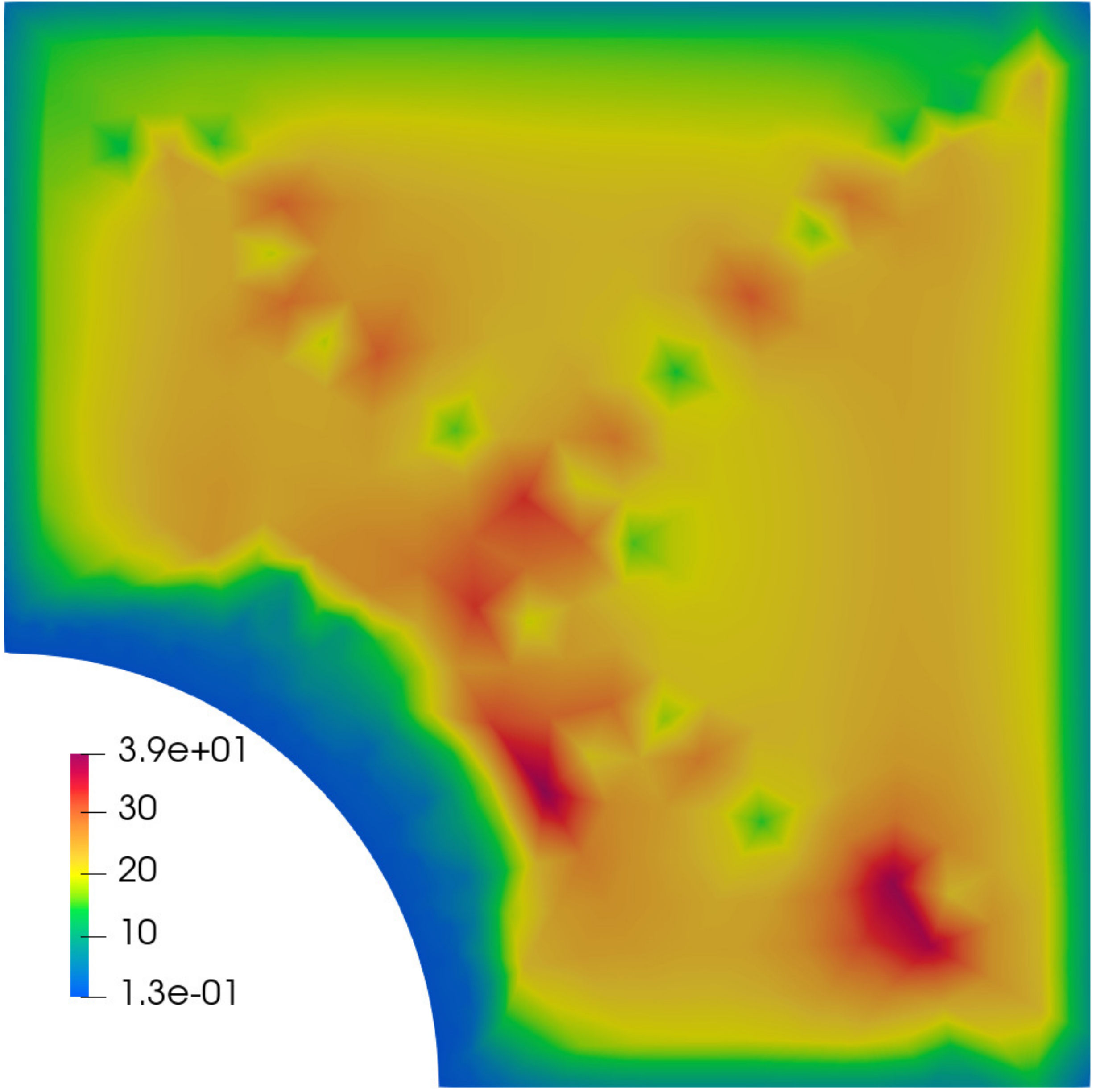} \\[\abovecaptionskip]
		\small (a) Consistent/discrete control sensitivities $dJ/d\boldsymbol{s}_\mathsmaller{\Gamma}$.
	\end{tabular}
	
	\vspace{\floatsep}
	
	\begin{tabular}{@{}c@{}}
		\includegraphics[width=0.7\linewidth,keepaspectratio]{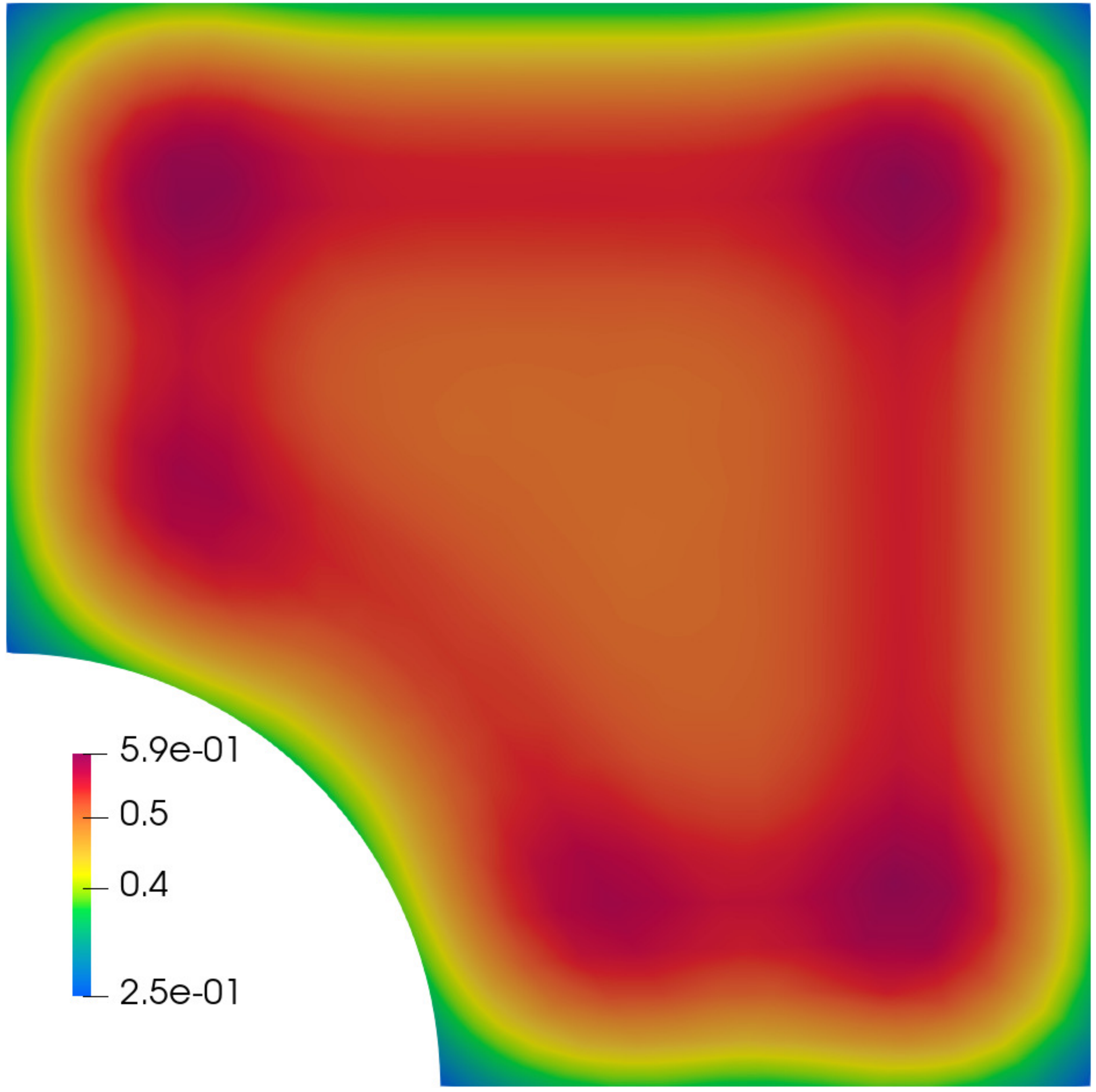} \\[\abovecaptionskip]
		\small (b) Scaled control sensitivities $dj/d\boldsymbol{s}_\mathsmaller{\Gamma} = \boldsymbol{M}_{\mathsmaller{\Gamma}}^{-1}  \,.\, dJ/d\boldsymbol{s}_\mathsmaller{\Gamma}$.
	\end{tabular}
	
	\caption{Control/filtered sensitivities computed with explicit filtering for the consistent/discrete sensitivities in Fig. \ref{fig:non-uni-plate}.}
	\label{fig:filtered-sens-explicit}
\end{figure}

We perform consistency tests with a uniform continuous sensitivity field of $dj/d\boldsymbol{x}_\mathsmaller{\Gamma}=[(0,0,1)]_{n_\mathsmaller{\Gamma}\times1}$ over a perforated plate which is discretized non-uniformly, Fig. \ref{fig:non-uni-plate}. 
It should be emphasized that although the continuous shape sensitivity field  $dj/d\boldsymbol{x}_\mathsmaller{\Gamma}$ is uniformly distributed, the consistent/discrete nodal sensitivities $dJ/d\boldsymbol{x}_\mathsmaller{\Gamma}$ are not necessarily uniform and indeed they are scaled by the mesh mass matrix, see Eq. \ref{e:discrete-continuous-sens}. Control sensitivities obtained after explicit filtering with a linear hat function are demonstrated in Fig. \ref{fig:filtered-sens-explicit}. As could be expected, the consistent/discrete control sensitivities are non-smooth and indeed mesh-dependent whereas the ones scaled by the inverse mass matrix are certainly smooth and mesh-independent. As a matter of fact, any consistent discrete field is mesh-dependent and noisy. Therefore, scaling the consistent sensitivities $dJ/d\boldsymbol{s}_\mathsmaller{\Gamma}$ with the inverse mass matrix is absolutely necessary to avoid mesh-dependency and introducing irrelevant noises into design. We also notice that the scaled sensitivities are non-uniform close to the boundary edges, indicating lack of consistency in mapping/filtering sensitivity with the explicit approach.

%

\section{Implicit shape filtering}
\label{s:Implicit shape filtering}
Filtering can also be performed implicitly by elliptic PDEs whose inverse operator is a local smoother. A well-established and commonly used implicit filter is the so-called Helmholtz/Laplace-Beltrami operator \citep{lazarov2011filters,kawamoto2011heaviside}, that is also labeled as the Sobolev smoothing \citep{Jameson2000,schmidt2008shape,mohammadi2009applied,dick2021combining}, which reads 

\begin{equation}
\label{e:sobolev-filter}
\begin{split}
\;\;\;\;\;\;& (I-\epsilon\,.\,\Delta)\, \widetilde{d} \,=\, d\\
\;\;\;&	\boldsymbol{n}.\nabla\widetilde{d} \,=\, 0, \;\;\; \text{along external boundaries.}
\end{split}
\end{equation}

where is $\widetilde{d}$ is the filtered field, $d$ is the raw field and $\epsilon$ is a scalar penalizing high spatial variations in $\widetilde{d}$. As mentioned previously, the so-called traction method \citep{azegami1994solution,azegami2006smoothing,riehl2014discrete} is also an implicit filtering scheme. This method determines the filtered field by solving a pseudo-elastic problem by the loading of the raw field. In this work, we develop surface and bulk-surface implicit smoothers based on the Helmholtz operator in Eq. \ref{e:sobolev-filter}. They are respectively applied to shell and solid geometries to ensure $C^1$ shape variations during the optimization iterations. It must be mentioned that the bulk-surface filtering treats boundary and internal domains simultaneously. Therefore, there is no need for the inclusion of any mesh motion technique in shape optimization for volumetric domains.      

\subsection{Surface filtering}
\label{s:Surface Helmholtz filtering}
The Sobolev/Helmholtz operator in Eq. \ref{e:sobolev-filter} can be casted on the design surface as 
\begin{equation}
	\label{e:belt-cont}
	-(r^{H}_\mathsmaller{\Gamma})^2\,\nabla_\mathsmaller{\Gamma}\,.\, \nabla_\mathsmaller{\Gamma} \,x^i+x^i  = s^i , \;\; i \in \{1,2,3\}; \;\; \text{on} \;\; \Gamma 
\end{equation}
$\nabla_\mathsmaller{\Gamma}\,.\, \nabla_\mathsmaller{\Gamma}$ is the Laplace-Beltrami operator. $\nabla_\mathsmaller{\Gamma}$ is the tangential gradient operator and can be calculated for a function $w$ defined on some open neighborhood of $\Gamma$ with a unit normal $\boldsymbol{n}$ as $\nabla_\mathsmaller{\Gamma} w = \nabla w -(\boldsymbol{n}. \nabla w)\, \boldsymbol{n} $, $r^{H}_\mathsmaller{\Gamma}$ is the filter radius/length parameter of Helmholtz-based filtering which plays a similar role as $r^{E}_\mathsmaller{\Gamma} $ in the explicit filtering and as it approaches zero, $r^{H}_\mathsmaller{\Gamma}\rightarrow0$, smoothing effect disappears, i.e. $\boldsymbol{x}=\boldsymbol{s}$. Smoothing property of Eq. \ref{e:belt-cont} in every spatial direction can be clearly explained by the minimization of the associated potential $\Pi$ as
\begin{equation}
	\Pi\left(x^i\right) = \frac{1}{2}\, (r^{H}_\mathsmaller{\Gamma})^2 \int_{\Gamma} |\nabla x^i|^2 \, d\Gamma + \frac{1}{2} \int_{\Gamma} \left(x^i-s^i\right)^2 \, d\Gamma
\end{equation} 
where the first and second integrals respectively measure the noisiness (spatial variation) of geometry and the difference between control and actual geometries. Here, $r^{H}_\mathsmaller{\Gamma}$ can be interpreted as the weight/importance of smoothness in the minimization of the  multi-objective function $\Pi$. It should be noted that Eq. \ref{e:belt-cont} is an isotropic in-surface filtering which does not depend on the direction. 

%

\subsection{Bulk-surface filtering}
\label{s:Bulk-surface Helmholtz filtering}
In node-based shape optimization with volumetric domains, we need to treat not only the surface/boundary nodes but also the internal nodes in order not to lose the mesh quality for numerical computations. In the literature so far, surface and internal nodes mainly have been treated separately in a master-follower manner, meaning that surface nodes are the master and internal nodes follow them using mesh motion techniques. The main issue with this approach is that not every smooth shape mode of the surface mesh can be tolerated by the volume mesh in terms of quality measures. Therefore, remeshing may be needed in the course of optimization process. To enforce the surface smoothness and maintain the volume mesh quality in one shot, we use an elliptic partial differential equation, posed on solid domain with smooth boundary, endowed with a generalized Robin boundary condition which involves the Laplace–Beltrami operator on the boundary surface as

\begin{subequations}
	\label{e:bulk_surface_PDE}
	\begin{align}
	-(r^{H}_\mathsmaller{\Omega})^2\;\nabla \;.\; \mathlarger{\boldsymbol{\sigma}}  + \boldsymbol{x} & = \boldsymbol{s}, \;\; && \text{in} \;\; \Omega \\
	-(r^{H}_\mathsmaller{\Gamma})^2\;\nabla_\mathsmaller{\Gamma}\,.\, \nabla_\mathsmaller{\Gamma} \boldsymbol{x} & = \boldsymbol{0} ,  \;\; && \text{on} \;\; \Gamma
	\end{align}
\end{subequations}

where $r^{H}_\mathsmaller{\Omega}$ is the Helmholtz bulk filter radius, $\nabla$ is the conventional spatial gradient operator, $\boldsymbol{\sigma}$ is a second order tensor similar to the Cauchy stress tensor in continuum mechanics and defined as 

\begin{subequations}
	\label{e:cauchy_bulk_surface_PDE}
	\begin{align}
&\mathlarger{\boldsymbol{\sigma}} =  \lambda \, \text{tr}(\mathlarger{\boldsymbol{\epsilon}}(\boldsymbol{x})) \, \boldsymbol{I} + 2\,\mu \, \mathlarger{\boldsymbol{\epsilon}}(\boldsymbol{x})\\
&\mathlarger{\boldsymbol{\epsilon}}(\boldsymbol{x}) = \frac{1}{2} \left(\nabla\boldsymbol{x}+\left(\nabla\boldsymbol{x}\right)^T\right)
	\end{align}
\end{subequations}
where $\text{tr}()$ is the trace operator, $\lambda$ and $\mu$ are the Lam\'e constants, $\boldsymbol{I}$ is the identity tensor, $\boldsymbol{\epsilon}$ is the strain tensor acting on geometry. Once again, one should note that $r^{H}_\mathsmaller{\Omega}$ and $r^{H}_\mathsmaller{\Gamma}$ are devised to control the surface and bulk filtering properties of the presented shape parameterization. To reveal these properties and their relation, we formulate Eqs. \ref{e:bulk_surface_PDE} as the minimization of associated potential 
\begin{equation}
\label{e:bulk-surface-energy}
\begin{aligned}
\Pi\left(\boldsymbol{x}\right) = \, & \frac{1}{2}\, (r^{H}_\mathsmaller{\Omega})^2 \int_{\Omega} \mathlarger{\boldsymbol{\sigma}} : \mathlarger{\boldsymbol{\epsilon}} \, d\Omega \,+\, \frac{1}{2} \, (r^{H}_\mathsmaller{\Gamma})^2 \int_{\Gamma} \nabla_\mathsmaller{\Gamma} \boldsymbol{x}^T . \nabla_\mathsmaller{\Gamma} \boldsymbol{x} \, d\Gamma \\ & + \frac{1}{2} \int_{\Omega} \left(\boldsymbol{x}-\boldsymbol{s}\right)^T .\left(\boldsymbol{x}-\boldsymbol{s}\right) \, d\Omega
\end{aligned}
\end{equation}
Here, we minimize a weighted sum of the total strain energy of the geometry occupying $\Omega$, the integrated spatial variation of the boundary geometry (noisiness) and the total deviation between the control and actual geometries. The corresponding weights are respectively $(r^{H}_\mathsmaller{\Omega})^2$, $(r^{H}_\mathsmaller{\Gamma})^2$ and 1, assuming uniform distribution over the integration domains. It should be noted that this is a multi-criteria optimization which may be posed to well-known challenges such as conflict between criteria and dominance over each other. Since smoothness of the design boundary $\Gamma$ (second term in Eq. \ref{e:bulk-surface-energy}) is of greater importance and to prevent it from being dominated by the volumetric strain (first term in Eq. \ref{e:bulk-surface-energy}), the Helmholtz bulk filter radius is calculated as follows

\begin{equation}
\label{e:bulk-surface-ratio}
	(r^{H}_\mathsmaller{\Omega})^2 = \beta \: \dfrac{(r^{H}_\mathsmaller{\Gamma})^2\mathlarger{\int_{\Gamma} \, \nabla_\mathsmaller{\Gamma} \boldsymbol{x}^T . \nabla_\mathsmaller{\Gamma} \boldsymbol{x} \, d\Gamma}}{\mathlarger{\int_{\Omega} \,\mathlarger{\boldsymbol{\sigma}} : \mathlarger{\boldsymbol{\epsilon}} \, d\Omega}}, \;\; 0 < \beta \leq 1 
\end{equation}

$\beta$ is a weighting factor which directly controls the importance of the volumetric strain minimization versus the boundary noisiness minimization. $\beta=1$ means that the volumetric and surface criteria are treated equally.    

\subsection{FEM discretization}
To solve PDEs, we use piecewise linear finite element functions to approximate the geometry as well as data fields, i.e. isoparametric finite elements. Therefore, actual geometry and its control field  within each element are approximated as $	x^i \approx \boldsymbol{N}^i \boldsymbol{x}^i,\;\; s^i \approx \boldsymbol{N}^i \boldsymbol{s}^i; \;\; i \in \{1,2,3\}$, where $\boldsymbol{N}^i$ is the vector of element shape functions in the $i$th Cartesian direction, $\boldsymbol{x}^i$ and $\boldsymbol{s}^i$ are the $i$th component of nodal coordinate vector and nodal control point vector, respectively. Furthermore, the traces of isoparametric bulk finite element functions on the boundary are used as surface finite elements. The same approach has been followed by \cite{edelmann2021isoparametric,dziuk2013finite} in the context of bulk–surface and surface PDEs. We derive a matrix–vector formulation of the discretized PDEs using a standard Galerkin-based finite element formulation. Then we can define following bulk and surface mass and stiffness matrices: 

\begin{subequations}
	\label{e:matrices}
	\begin{align}
	&\boldsymbol{M}_\mathsmaller{\Omega} = \sum\int_{\Omega^e} \boldsymbol{N}^T\boldsymbol{N} d\Omega,\\ 		&\boldsymbol{M}_\mathsmaller{\Gamma} = \sum\int_{\Gamma^e} \boldsymbol{N}^T_\mathsmaller{\Gamma}\boldsymbol{N}_\mathsmaller{\Gamma} d\Gamma \\
	&\boldsymbol{K}_\mathsmaller{\Omega} = \sum\int_{\Omega^e} (r^{H}_\mathsmaller{\Omega^e})^2\boldsymbol{B}^T\boldsymbol{C}\boldsymbol{B} \,d\Omega \\ &\boldsymbol{K}_\mathsmaller{\Gamma} = (r^{H}_\mathsmaller{\Gamma})^2\sum\int_{\Gamma^e} (\nabla_\mathsmaller{\Gamma}\boldsymbol{N}_\mathsmaller{\Gamma})^T\nabla_\mathsmaller{\Gamma}\boldsymbol{N}_\mathsmaller{\Gamma} d\Gamma				
	\end{align}
\end{subequations}
where $\boldsymbol{B}$ contains the spatial gradients of the bulk shape functions  $\boldsymbol{N}$, $\boldsymbol{C}$ is the linear-elastic isotropic constitutive matrix, $\boldsymbol{N}_\mathsmaller{\Gamma}$ is the traces of bulk shape functions on the boundary, $r^{H}_\mathsmaller{\Omega^e}$ is the elemental Helmholtz bulk filter radius which is assumed to be spatially varying. In this work, with the aim to maintain the mesh quality and reduce the frequency of remeshing during the shape optimization, the bulk filter radius is selected element-wise, based on element sizes. Among others, “Mesh-Jacobian-based stiffening” (MJBS) \citep{tezduyar,stein} has been widely used to selectively stiffen or soften elements against shape changes. It can be implemented simply by dropping Jacobian from the finite element formulation of the mesh governing equations, resulting in the smaller elements being stiffened more than the larger ones. That can be realized by choosing the elemental Helmholtz filter radius as $\mathlarger{r^{H}_\mathsmaller{\Omega^e} = \frac{J^0}{J^e}}$, where $J^e$ is the Jacobian for element $e$ and $J^0$ is a scaling parameter which is calculated based on Eq. \ref{e:bulk-surface-ratio} as follows

\begin{equation}
J^0 = \beta \: \dfrac{(r^{H}_\mathsmaller{\Gamma})^2 \, (\boldsymbol{x}_\mathsmaller{\Gamma})^T \, \boldsymbol{K}_\mathsmaller{\Gamma} \, \boldsymbol{x}_\mathsmaller{\Gamma} }{(\boldsymbol{x})^T \, \mathlarger{\left(\sum\int_{\Omega^e} (\frac{1}{J^e})\,\boldsymbol{B}^T\boldsymbol{C}\boldsymbol{B} \,d\Omega\right)} \: \boldsymbol{x}}
\end{equation}

Having evaluated the bulk and surface mass and stiffness matrices, the shape governing equations for shells and solids (Eqs. \ref{e:belt-cont},\ref{e:bulk_surface_PDE}) read respectively as
\begin{subequations}
	\label{e:matrix-vector-format}
	\begin{align}
	\boldsymbol{x}_\mathsmaller{\Gamma} & = (\boldsymbol{A}^{I}_\mathsmaller{\Gamma})^{-1} \,.\, \boldsymbol{s}_\mathsmaller{\Gamma}; && \boldsymbol{A}^{I}_\mathsmaller{\Gamma} = (\boldsymbol{K}_\mathsmaller{\Gamma}\,+\,\boldsymbol{M}_\mathsmaller{\Gamma})^{-1}\,.\,\boldsymbol{M}_\mathsmaller{\Gamma}
	\\
	\boldsymbol{x} &= (\boldsymbol{A}^{I})^{-1}   \,.\, \boldsymbol{s}; && \boldsymbol{A}^{I}_\mathsmaller{\Gamma} = (\boldsymbol{K}_\mathsmaller{\Gamma}\,+\,\boldsymbol{K}_\mathsmaller{\Omega}\,+\,\boldsymbol{M}_\mathsmaller{\Omega})^{-1} \,.\,\boldsymbol{M}_\mathsmaller{\Omega}
	\end{align}
\end{subequations}

where the superscript $I$ refers to the implicit filtering,  $\boldsymbol{A}^{I}_\mathsmaller{\Gamma}$ and $\boldsymbol{A}^{I}$ are the implicit filter matrix of surface and bulk geometries, respectively. As the corresponding tangent matrices are symmetric and positive definite, the solution of the linear system can be achieved by utilizing efficient iterative solvers such as Conjugate Gradient.

\subsection{Sensitivity analysis and update rule}
Straight forward, applying the chain rule of differentiation the derivative of the response function with respect to the surface and volumetric control points is given respectively as
\begin{subequations}
	\label{e:implicit-sensitivity-mapping}
	\begin{align}
	\frac{dJ}{d\boldsymbol{s}_\mathsmaller{\Gamma}} & = \boldsymbol{M}_\mathsmaller{\Gamma} \,.\, (\boldsymbol{K}_\mathsmaller{\Gamma}\,+\,\boldsymbol{M}_\mathsmaller{\Gamma})^{-1}\,.\, \frac{dJ}{d\boldsymbol{x}_\mathsmaller{\Gamma}}\\
	\frac{dJ}{d\boldsymbol{s}} & = \boldsymbol{M}_\mathsmaller{\Omega} \,.\, (\boldsymbol{K}_\mathsmaller{\Gamma}\,+\,\boldsymbol{K}_\mathsmaller{\Omega}\,+\,\boldsymbol{M}_\mathsmaller{\Omega})^{-1} \,.\, \frac{dJ}{d\boldsymbol{x}}
	\end{align}
\end{subequations}

It should be reminded that discrete sensitivities in the above equations are consistent nodal values which show size effects in the gradients. To avoid mesh dependency, the discrete control sensitivities are scaled by the inverse mass matrix, which basically means that we reconstruct the control gradient field sampled at mesh points from the consistent nodal sensitivities, viz. 

\begin{subequations}
	\label{e:implicit-scaled-control-sens}
	\begin{align}
	\frac{dj}{d\boldsymbol{s}_\mathsmaller{\Gamma}} & = \boldsymbol{M}^{-1}_\mathsmaller{\Gamma} \,.\,\frac{dJ}{d\boldsymbol{s}_\mathsmaller{\Gamma}}=(\boldsymbol{K}_\mathsmaller{\Gamma}\,+\,\boldsymbol{M}_\mathsmaller{\Gamma})^{-1}\,.\, \frac{dJ}{d\boldsymbol{x}_\mathsmaller{\Gamma}}\\
	\frac{dj}{d\boldsymbol{s}} & = \boldsymbol{M}^{-1}_\mathsmaller{\Omega} \,.\,\frac{dJ}{d\boldsymbol{s}}= (\boldsymbol{K}_\mathsmaller{\Gamma}\,+\,\boldsymbol{K}_\mathsmaller{\Omega}\,+\,\boldsymbol{M}_\mathsmaller{\Omega})^{-1} \,.\, \frac{dJ}{d\boldsymbol{x}}
	\end{align}
\end{subequations}

The scaled gradients do not include disturbing discretization effects and ensure computations of efficient search directions without the necessity of second order information. Therefore, the discrete control field and subsequently the geometry are updated as follows
\begin{equation}
\label{e:quasi-newton-update-rule-implicit}
\begin{split}
\Delta \boldsymbol{s}_\mathsmaller{\Gamma} \, & =\, -\alpha  \,.\, \frac{dj}{d\boldsymbol{s}_\mathsmaller{\Gamma}} && \rightarrow \Delta \boldsymbol{x}_\mathsmaller{\Gamma} = \boldsymbol{A}^{I}_\mathsmaller{\Gamma} \,.\, \Delta \boldsymbol{s}_\mathsmaller{\Gamma}\\
\Delta \boldsymbol{s} \, & =\, -\alpha  \,.\, \frac{dj}{d\boldsymbol{s}} && \rightarrow \Delta \boldsymbol{x} = \boldsymbol{A}^{I} \,.\, \Delta \boldsymbol{s}
\end{split}
\end{equation}

As mentioned earlier, one can interpret the proposed control update rule using the scaled gradients as a quasi-Newton step with the diagonal approximation of the Hessian matrix by the mass matrix.

\subsection{Consistency check}
\begin{figure}
	\centering
	\begin{tabular}{@{}c@{}}
		\includegraphics[width=0.7\linewidth,keepaspectratio]{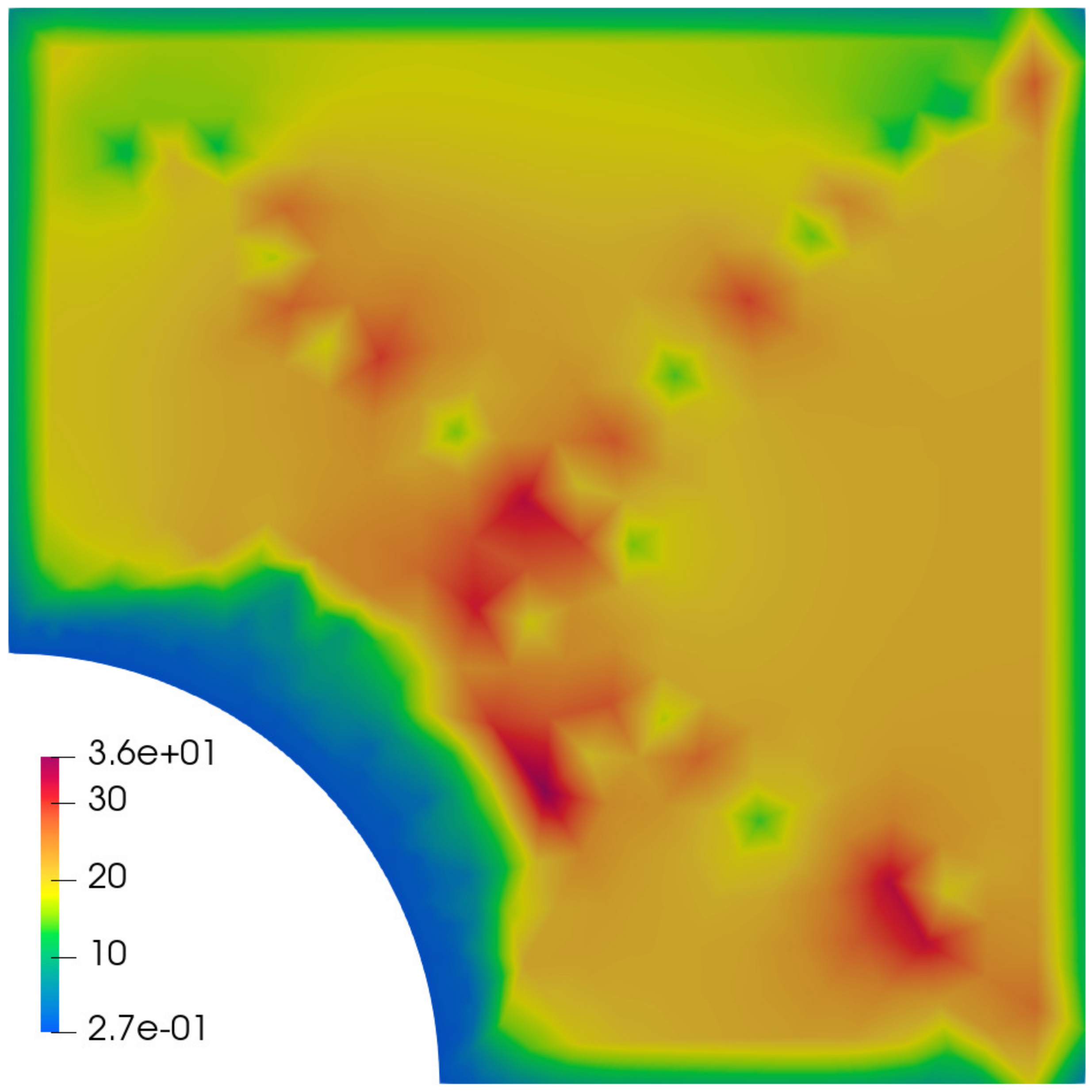} \\[\abovecaptionskip]
		\small (a) Discrete control sensitivities $dJ/d\boldsymbol{s}_\mathsmaller{\Gamma}$.
	\end{tabular}
	
	\vspace{\floatsep}
	
	\begin{tabular}{@{}c@{}}
		\includegraphics[width=0.7\linewidth,keepaspectratio]{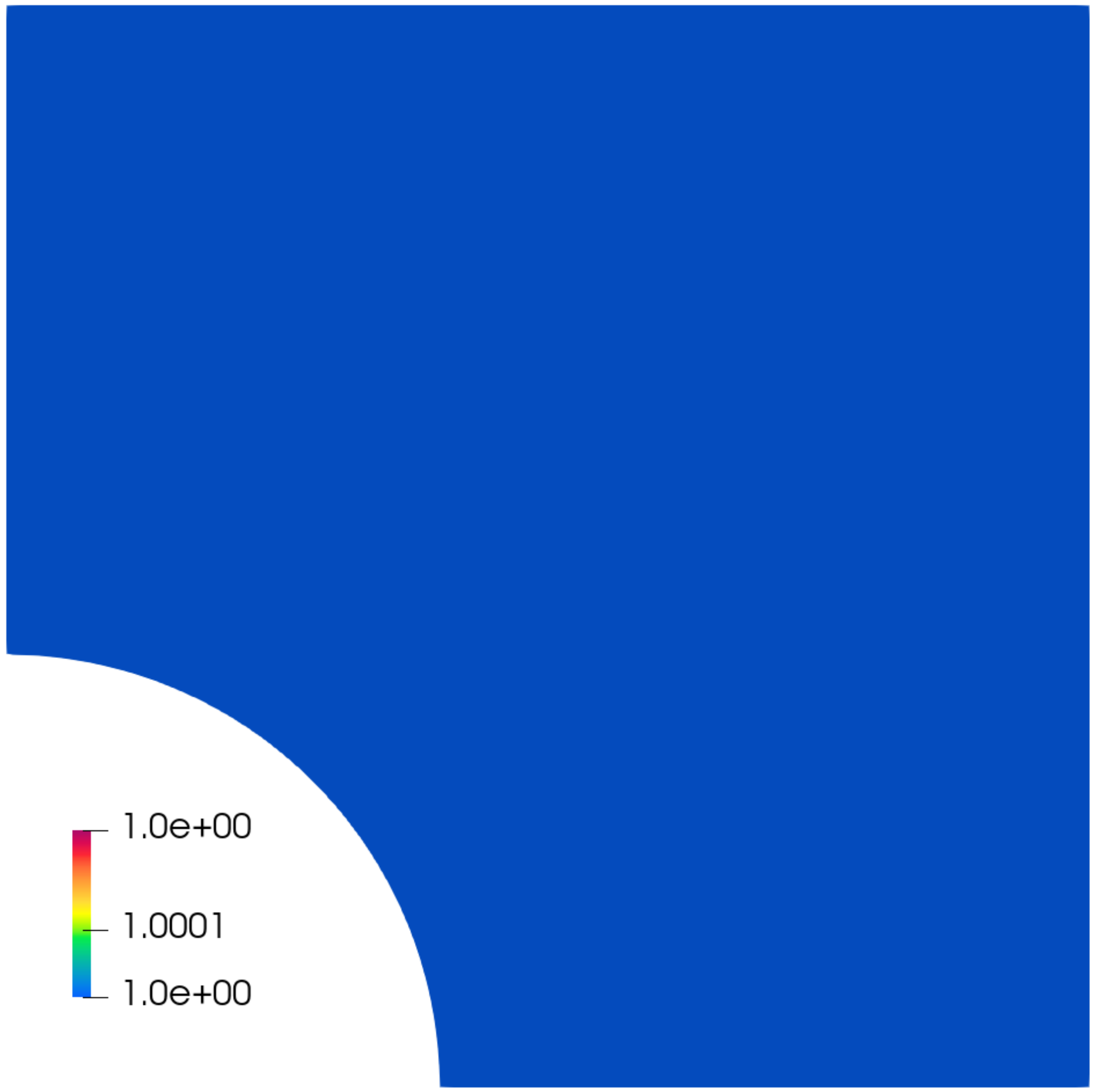} \\[\abovecaptionskip]
		\small (b) Scaled control sensitivities $dj/d\boldsymbol{s}_\mathsmaller{\Gamma} = \boldsymbol{M}_{\mathsmaller{\Gamma}}^{-1}  \,.\, dJ/d\boldsymbol{s}_\mathsmaller{\Gamma}$.
	\end{tabular}
	
	\caption{Control/filtered sensitivities computed with implicit filtering for the consistent/discrete sensitivities in Fig. \ref{fig:non-uni-plate}.}
	\label{fig:filtered-sens-implicit}
\end{figure}

The discrete implicit filtering was derived and discussed using conventional finite elements. As is well known in finite element theory, a proper finite element formulation includes the rigidbody motion capability. However, for the sake of completeness, we numerically study consistency aspects of the implicit filtering. Perforated plate used in Section \ref{s:explicit-consistency-check} is reconsidered and the synthetic consistent sensitivity field in Fig. \ref{fig:non-uni-plate} is used for the test. Control sensitivities obtained after implicit (Sobolev/Helmholtz) filtering are demonstrated in Fig. \ref{fig:filtered-sens-implicit}. Similar to the explicit filtering, the consistent/discrete control sensitivities are highly mesh-dependent whereas the scaled ones are fully uniform and mesh-independent. The obtained results indicate that the consistency (rigid-body-movement production) can be perfectly achieved throughout the implicit shape filtering process if control sensitivities are scaled by the inverse mesh mass matrix, independent of the discretization.

\section{Comparison of the explicit and implicit filters}
Relation between the convolution-based (explicit) and the Sobolev-/Helmholtz-based (implicit) filters has been discussed in the literature. \cite{stuck2011adjoint} state that an explicit filtering using a Gaussian kernel function is a first-order approximation to the Sobolev-based filtering. On the other hand, \cite{lazarov2011filters} have shown in 1D the correlation between the explicit filtering using linear hat function and the Helmholtz-based filtering. The fundamental solution of the non-homogeneous modified Helmholtz equation in Eq. \ref{e:belt-cont} for an infinite domain reads \citep{polyanin2016handbook}

\begin{equation}
\label{e:non-reg_Green}
\begin{aligned}
& \boldsymbol{x}_{0} = \int\limits_{-\infty}^{\infty}\int\limits_{-\infty}^{\infty}  \int\limits_{-\infty}^{\infty} F^H(\boldsymbol{ x},\boldsymbol{ x}_{0}) \: \boldsymbol{s}(\boldsymbol{ x}) \: dx^1\, dx^2 \,dx^3 \\
& F^H(\boldsymbol{ x},\boldsymbol{ x}_{0}) = \frac{1}{4\pi\norm{\boldsymbol{ x}-\boldsymbol{ x}_{0}}/(r^H)^2} e^{-\norm{\boldsymbol{ x}-\boldsymbol{ x}_{0}}/r^H}
\end{aligned}
\end{equation}

which is basically a convolutional filtering using Green's function $F^H$ as the kernel. It should be mentioned that the integral of Green's function is always one, therefore the consistency property is naturally included. We also note that the Green's function in Eq. \ref{e:non-reg_Green} is singular at the evaluation point $\boldsymbol{ x}_{0}$, therefore we propose a regularized  form of it as 
\begin{equation}
\label{e:reg_Green}
F^{RH}(\boldsymbol{ x},\boldsymbol{ x}_{0}) = \frac{1}{1+(4\pi\norm{\boldsymbol{ x}-\boldsymbol{ x}_{0}}/(r^H)^2)} e^{-\norm{\boldsymbol{ x}-\boldsymbol{ x}_{0}}/r^H}
\end{equation}

\begin{figure}
	\centering
	\includegraphics[keepaspectratio,width=0.5\textwidth]{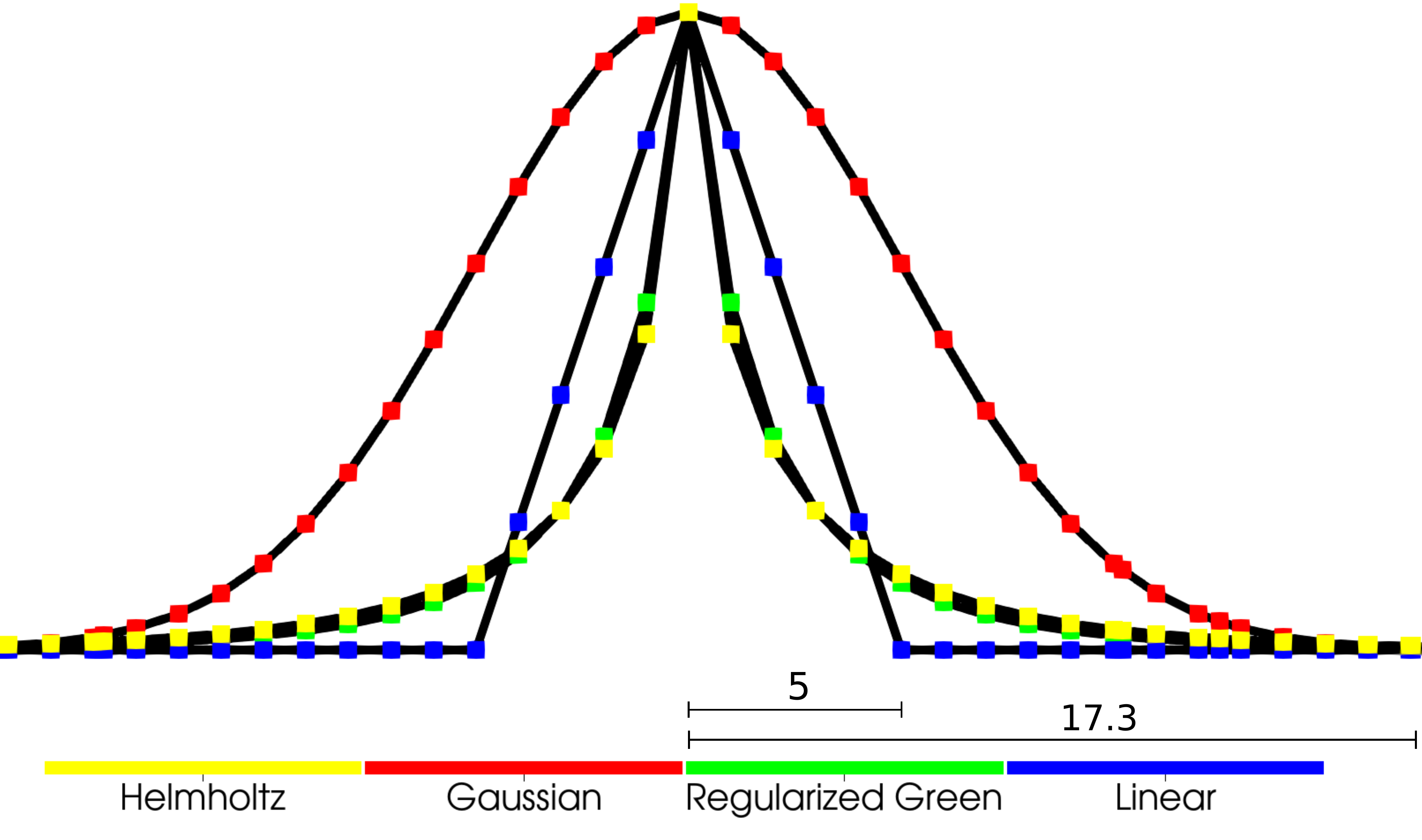}
	\caption{Discretized profiles of the considered filter functions.}
	\label{fig:kernels_comp}
\end{figure}

In order to compare qualitatively explicit and implicit filters, a square flat plate of $100\times 100$ meshed by triangles with the average element size of 1 is used. Figure \ref{fig:kernels_comp} illustrates schematically the considered filter functions for a given radius of 5. It should be mentioned that the Helmholtz's kernel at every mesh point is calculated numerically from the corresponding row in the inverted Helmholtz's stiffness matrix, i.e. $(\boldsymbol{K}_\mathsmaller{\Gamma}\,+\,\boldsymbol{M}_\mathsmaller{\Gamma})^{-1}$. The first observation indicates that the proposed regularized Green's function almost matches with the numerical kernel of the Helmholtz filter. We also notice that the linear hat function is approximately spanning over  $1/2\sqrt{3}$ of the other's support span, which confirms the relation developed by \cite{lazarov2011filters}. Furthermore, during the upcoming sections, performance of the presented filters is critically discussed for shape optimization of shells and solids receptively.

\subsubsection{Shell optimization}
\label{sec:shell_opt}
\begin{figure}
	\centering
	\includegraphics[keepaspectratio,width=0.5\textwidth]{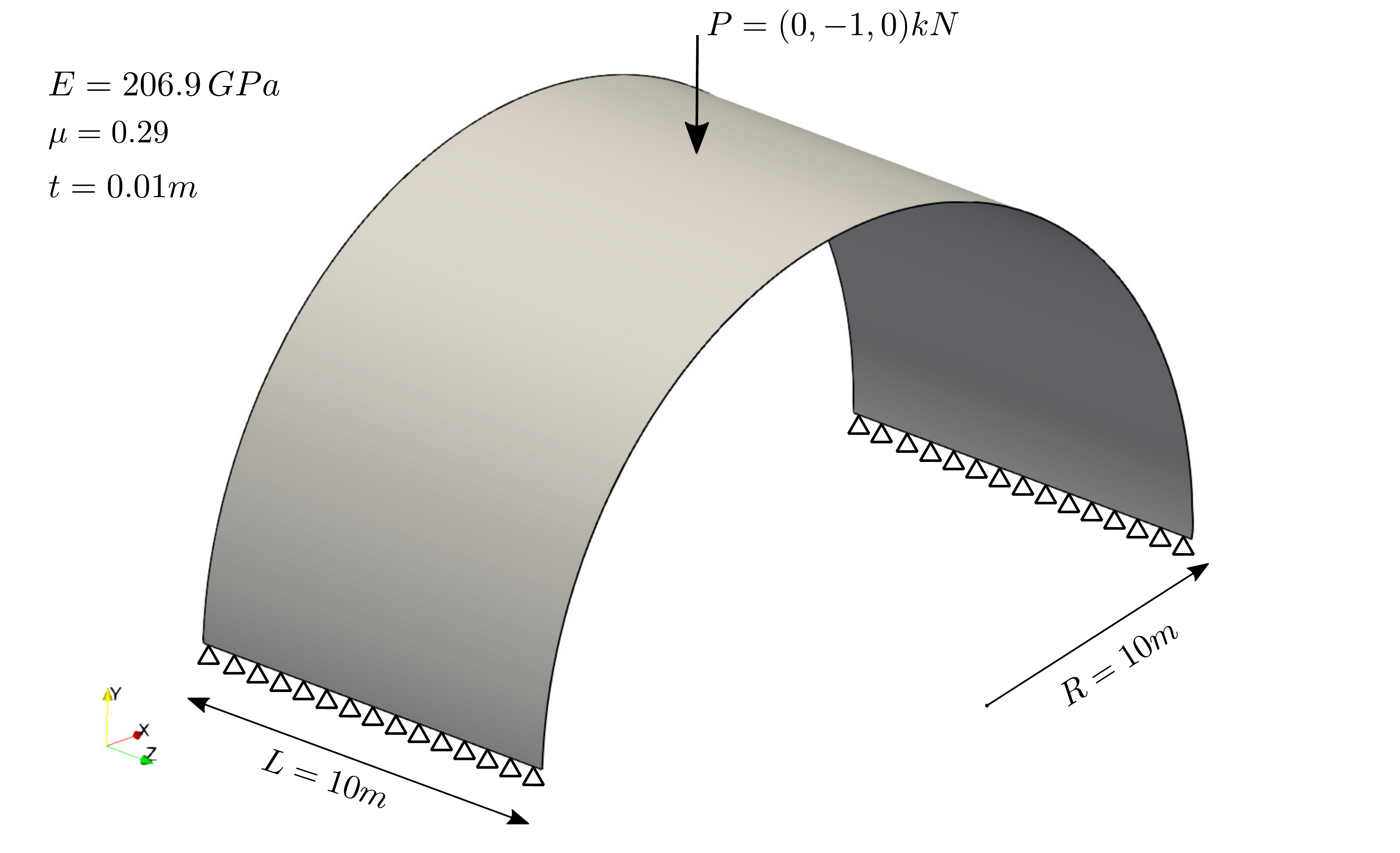}
	\caption{Geometry, support and loading of half-cylinder shell.}
	\label{fig:half-cylinder-problem-def}
\end{figure}

This example intends to investigate the influence of filter kernels on node-based shape optimization, in terms of numerical properties and the optimal shape. Fig. \ref{fig:half-cylinder-problem-def} shows a half-cylinder shell fully clamped at edges and loaded by a nodal force at the center. The specified geometry is described by 73522 unstructured triangles with the average element size of $0.1m$. We have implemented this study with the \emph{ShapeOptimization} and \emph{StructuralMechanics} applications of the open-source software \emph{Kratos-Multiphysics} \citep{dadvand2013migration}. As a linear system solver, we use AMGCL - an efficient parallel iterative linear solver \cite{Demidov2020}. 
 
\begin{figure}
	\centering
	\begin{tabular}{@{}c@{}}
		\includegraphics[width=\linewidth,keepaspectratio]{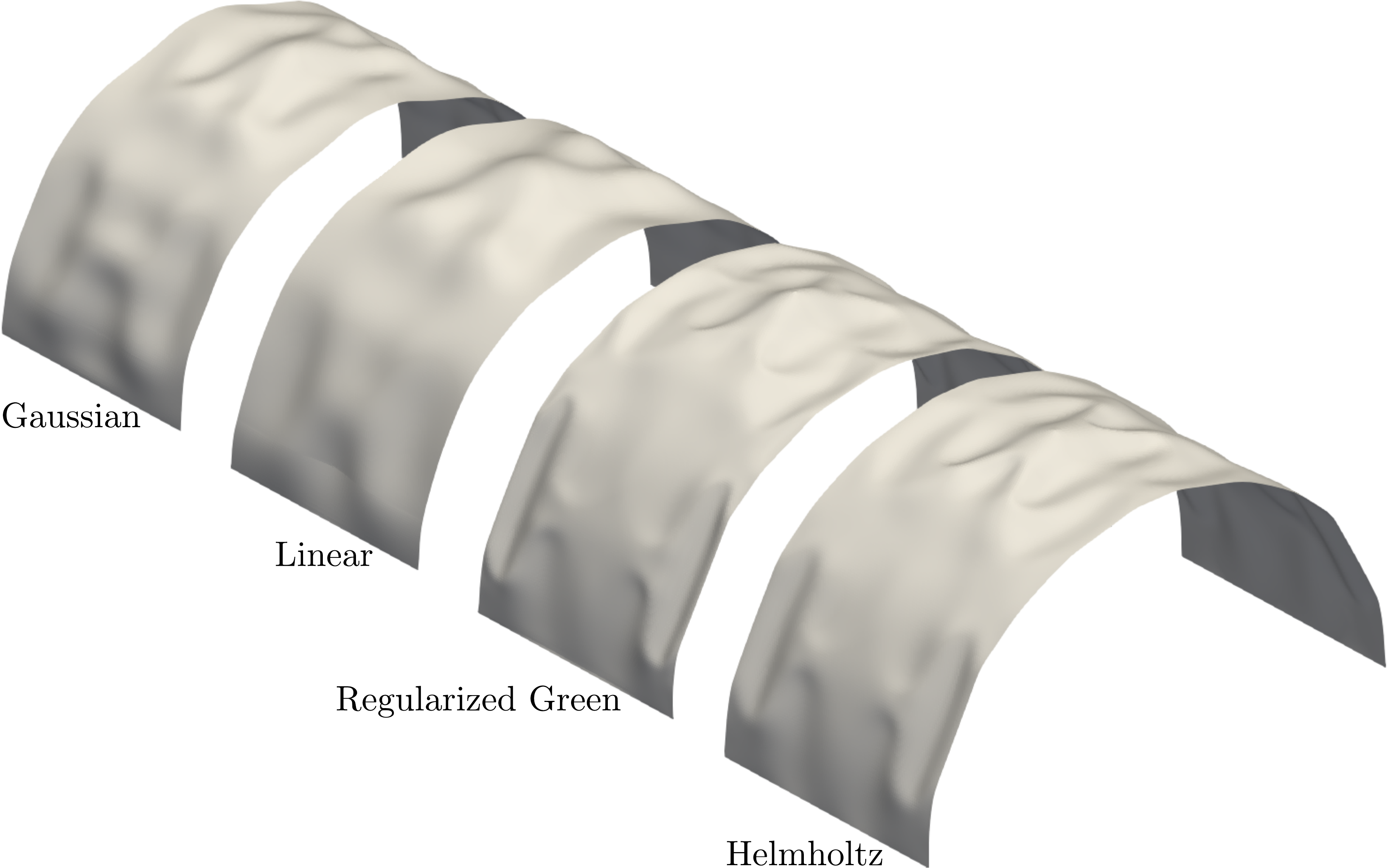} \\[\abovecaptionskip]
		\small (a) Optimal geometries obtained with different kernels\\ of the same support span $p=2m$.
	\end{tabular}

	\vspace{\floatsep}

	\begin{tabular}{@{}c@{}}
		\includegraphics[width=\linewidth,keepaspectratio]{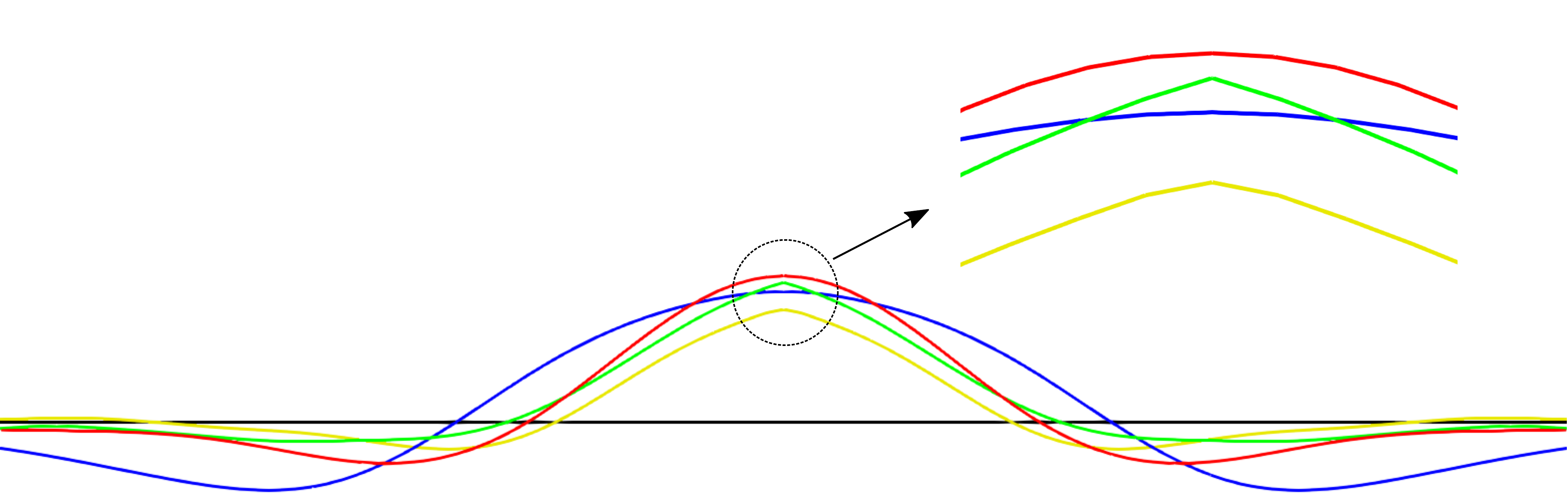} \\[\abovecaptionskip]
		\small (b) Cross sections of optimal geometries through the \\ YZ plane at the center.
	\end{tabular}

	\vspace{\floatsep}

	\begin{tabular}{@{}c@{}}
		\includegraphics[width=0.8\linewidth,keepaspectratio]{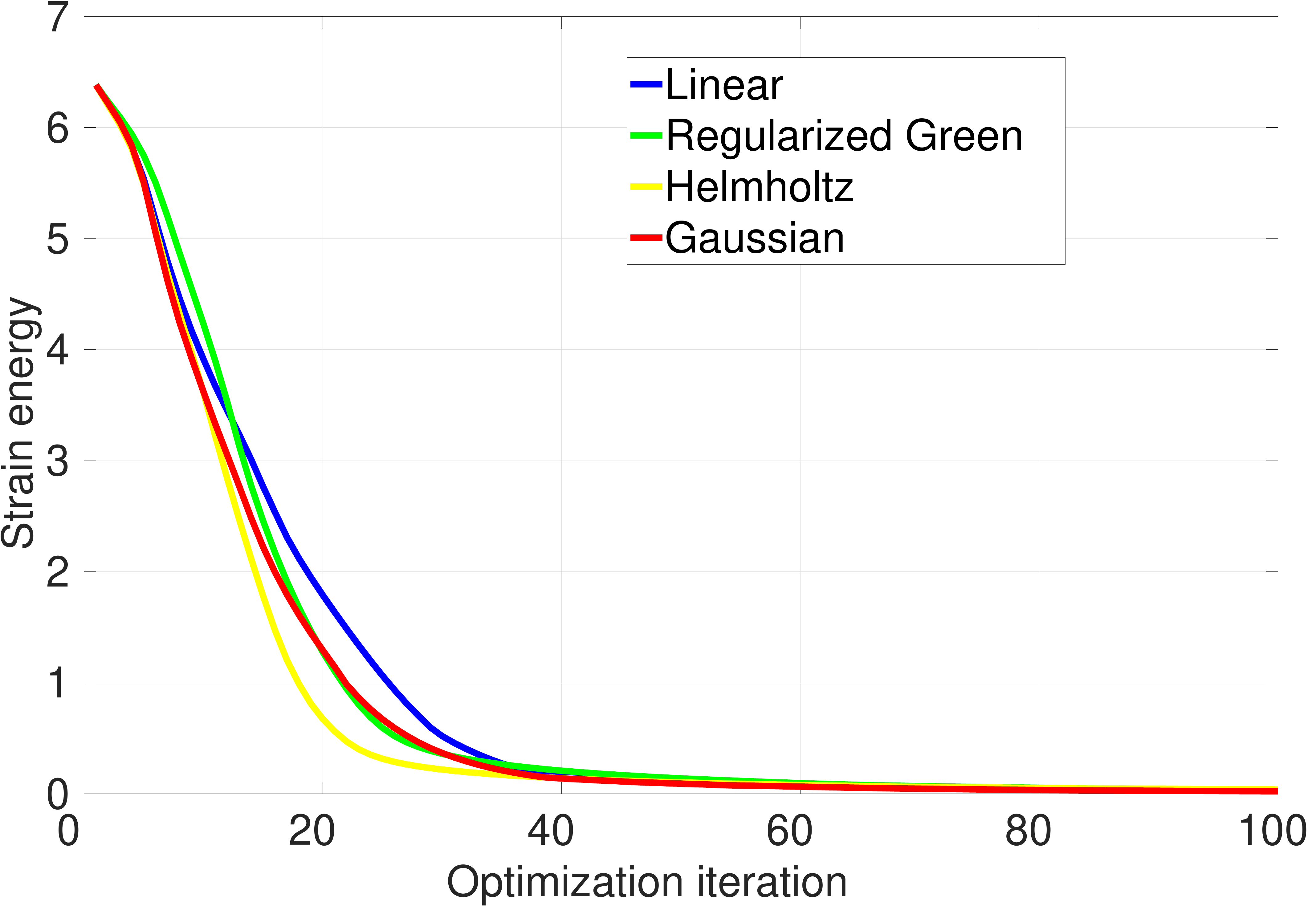} \\[\abovecaptionskip]
		\small (c) Convergence of objective.
	\end{tabular}

	\caption{Optimal design of the half-cylinder shell.}
	\label{fig:optimal-design-shell}
\end{figure}

Helmholtz, regularized Green, Gaussian and linear hat kernel functions are used for the following investigations. The projected steepest descent algorithm with constant step size \citep{asl2019shape} is employed to minimize the linear strain energy of the structure while constraining its mass. Optimal geometries obtained for the considered kernel function are presented in Fig.\ref{fig:optimal-design-shell}a. Filter radii are set so that the kernel functions have the same support span $p=2m$. It is observed that the Helmholtz and regularized Green kernels resulted in similar optimal geometries which are sharper than the others due to the resolution of local features. On the other hand, it can be seen from Fig.\ref{fig:optimal-design-shell}b that, the Helmholtz and regularized Green kernels allow the creation of a kink at the center; therefore not only do they obtain smooth shapes but could develop a kink as well \citep{muller2021novel}. From these observations it is reasonable to conclude that the Gaussian and linear kernel functions act as a low-pass filter suppressing fine scales, whereas the Helmholtz and regularized Green kernels are high-pass smoother which allows generation of fine scales. The improved objective function is depicted in Fig.\ref{fig:optimal-design-shell}c. It shows that all filter functions reduced the structural strain energy significantly, however optimization with the Helmholtz filtering seems to convergence faster.

\begin{figure}
	\centering
	\begin{tabular}{@{}c@{}}
		\includegraphics[width=\linewidth,keepaspectratio]{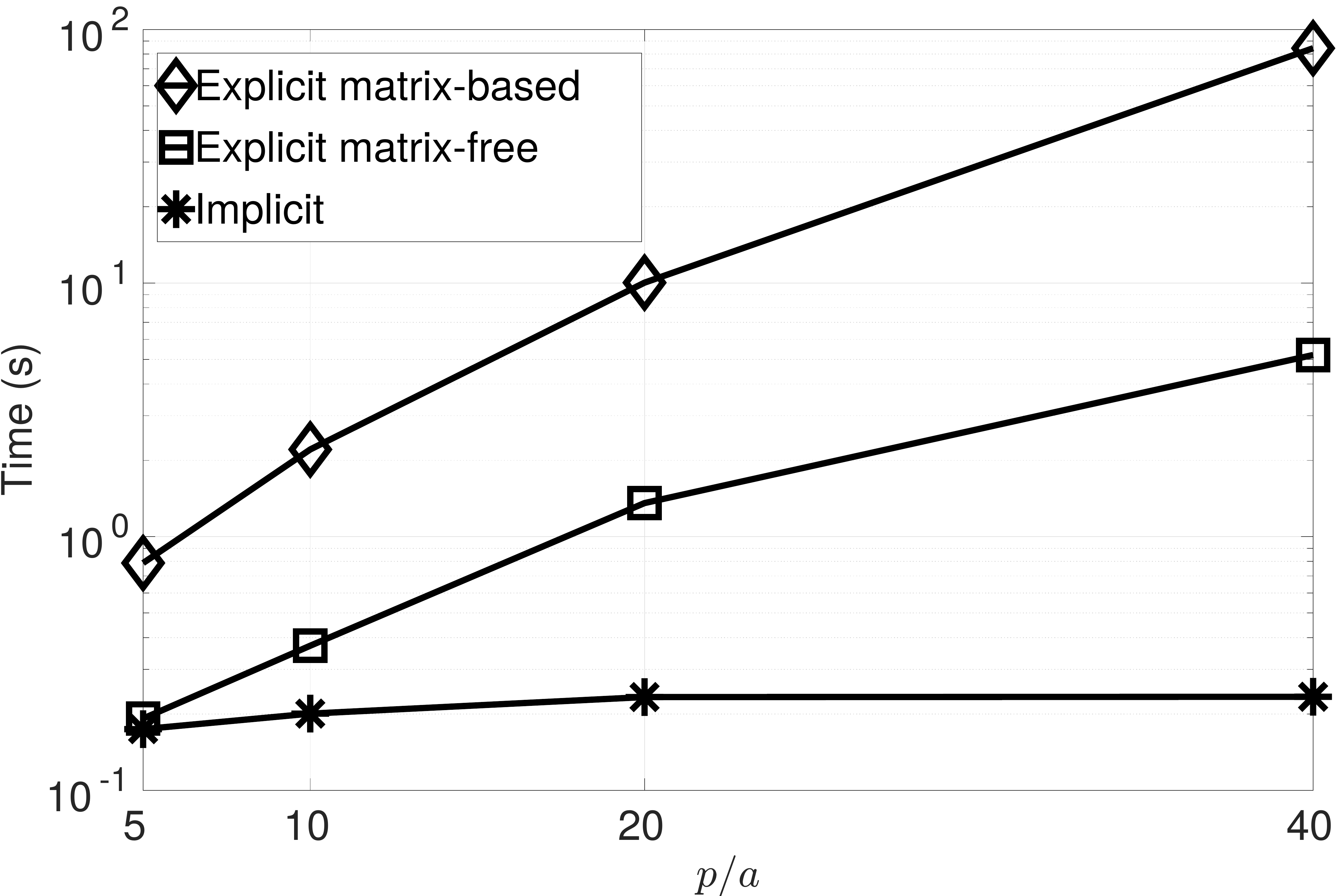} \\[\abovecaptionskip]
		\small (a) Time necessary for computing filtered field using \\implicit (Helmholtz) filter and the explicit filter with \\ (matrix-based) and without (matrix-free) storing the filter \\matrix for different support span to element size ratios $p/a$.
		
	\end{tabular}

	\vspace{\floatsep}
	
	\begin{tabular}{@{}c@{}}
		\includegraphics[width=\linewidth,keepaspectratio]{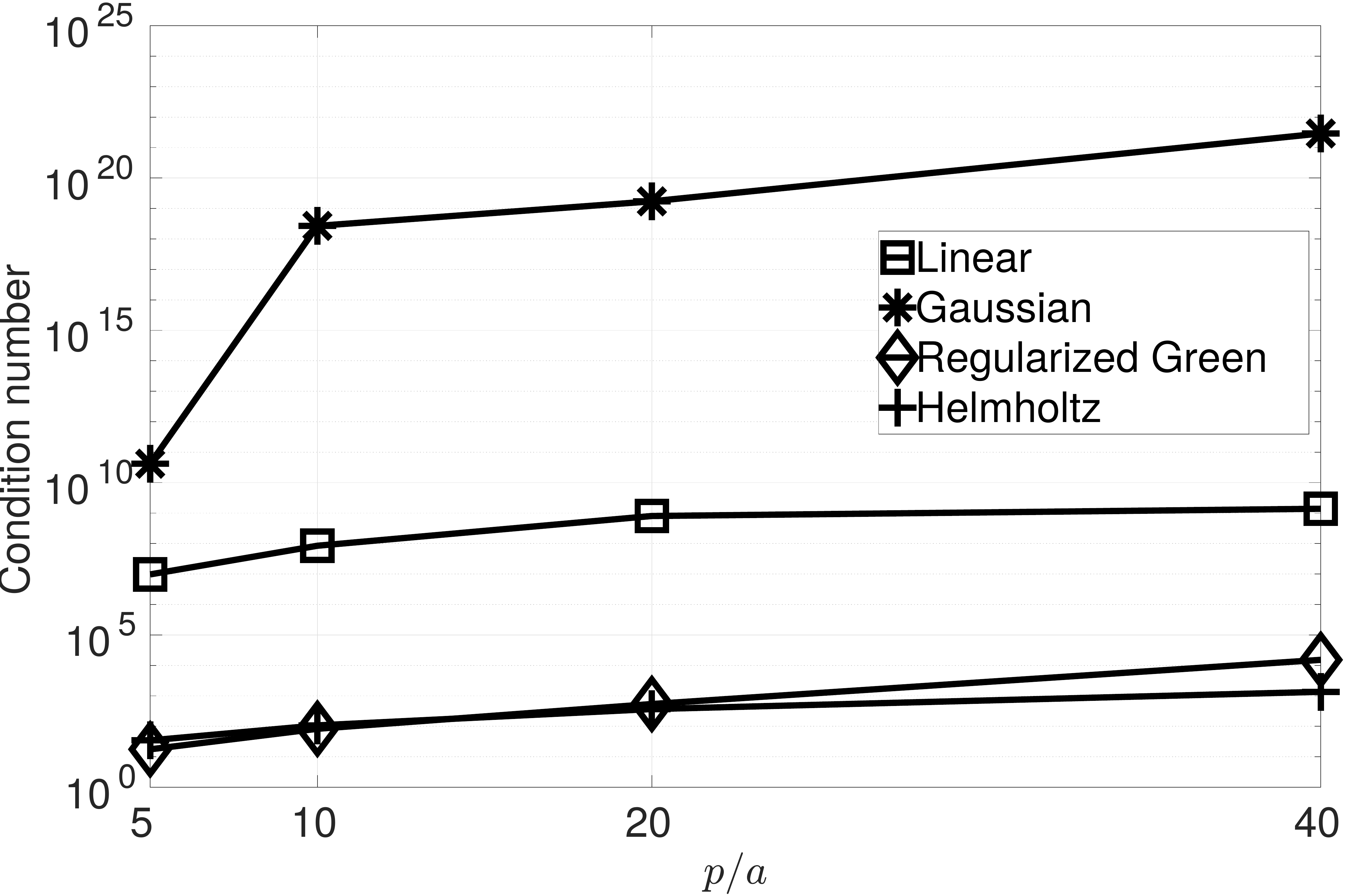} \\[\abovecaptionskip]
		\small (b) Condition number of the discrete filter operator\\ matrix of the considered kernels for different support \\ span to element size ratios $p/a$.
	\end{tabular}
	
	\caption{Numerical properties of different kernel functions.}
	\label{fig:numerical-properties-shells}
\end{figure}

Comparison of the time necessary for applying the explicit and implicit filters for different support span to element size ratios $p/a$ is shown in Fig. \ref{fig:numerical-properties-shells}a. For reference, this study is performed on a 6 core XEON W-2133 3.60GHz processors using hyper-threads-based parallelism. The neighbor search and distance evaluation needed for the explicit filtering is done by use of octrees for the sake of efficiency. Implementation of the explicit filter with (matrix-based) and without (matrix-free) storing the filter matrix is slower than the implicit PDE filter. Cost of the explicit filters grows almost linearly with the $p/a$ ratio, whereas that of the implicit filter is not growing dramatically and for high ratios it is almost constant. In accordance with \cite{lazarov2011filters}, overall conclusion is that the cost of building a search tree and performing neighbor search to explicity filter a surface field is more than the cost of a linear system solve associated with the implicit surface filtering. 

To investigate the sensitivity of filtering and subsequently the optimal shape to the filter radius, condition number of the discrete filter operator is used in this work. This is particularly important since the conditioning of the filter matrix directly influences the conditioning of the optimization problem, see \cite{najian2017consistent} for discussions and studies on this subject. Conditioning of the discrete filter operator matrix is shown in Fig. \ref{fig:numerical-properties-shells}b as a function of the support span. It was predictable that the Helmholtz and regularized Green kernel functions would behave very similarly. They show very good-conditioning with much less sensitivity to the $p/a$ ratio. Whereas, the linear and Gaussian kernels have generally ill-conditioning behavior which gets worse with increasing the support span. So we can conclude that the regularized Green kernel function may be preferred over the others for explicit shape filtering.   	  

\subsubsection{Solid Optimization}   
\begin{figure}
	\centering
	\begin{tabular}{@{}c@{}}
		\includegraphics[width=\linewidth,keepaspectratio]{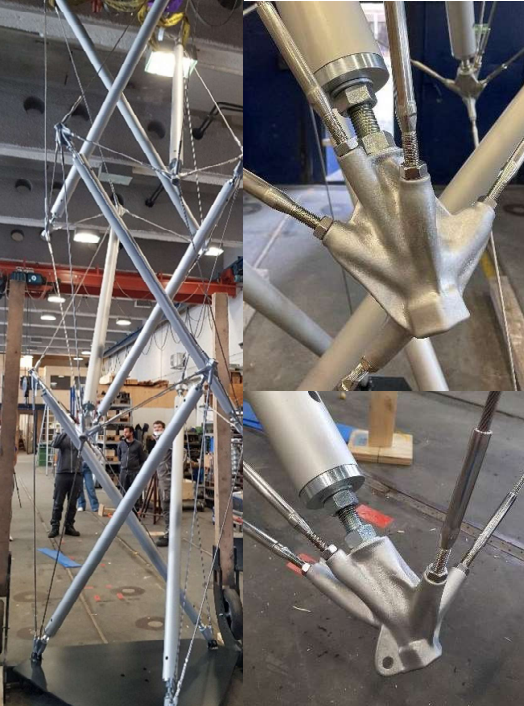} \\[\abovecaptionskip]
		\small (a) The five-meter “tensegrity tower” exhibited in the Deutsches \\ Museum of Science and Technology in Munich.
	\end{tabular}

	\vspace{\floatsep}
	
	\begin{tabular}{@{}c@{}}
		\includegraphics[width=\linewidth,keepaspectratio]{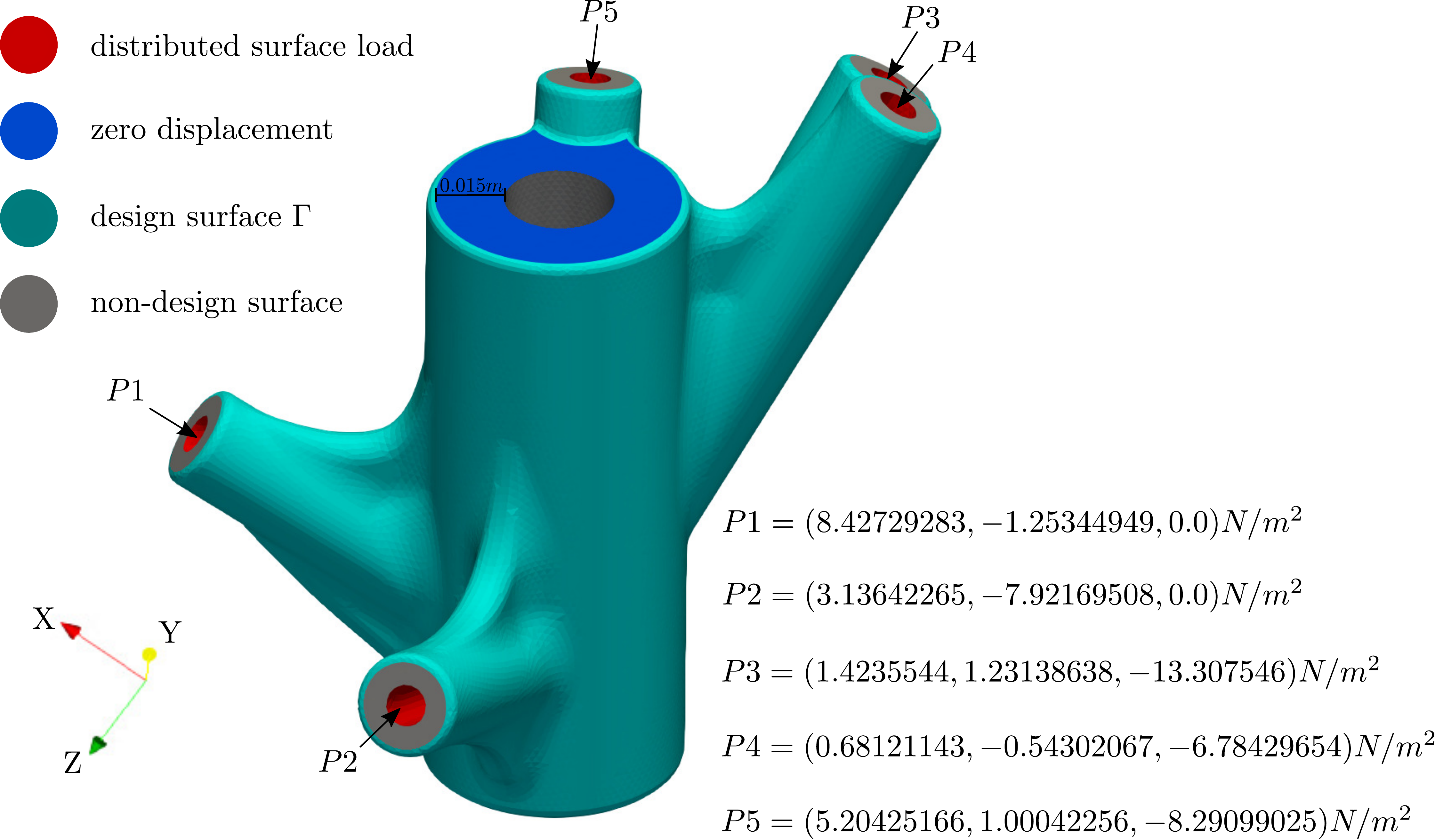} \\[\abovecaptionskip]
		\small (b) Problem setup of an aluminum connecting node.
	\end{tabular}

	\caption{Shape optimization of the nodes of a tensegrity-tower.}
	\label{fig:tensegrity-tower}
\end{figure}

In \cite{ghantasala2021node}, the Vertex-Morphing technique as an explicit surface filtering was used in shape optimization of the connecting nodes of a five-meter tensegrity (“tensional integrity”) tower, see Fig. \ref{fig:tensegrity-tower}a. The tower is float in the air and built out of tensile and compressive elements. The nodes are highly complex connections between the rods and cables and they are additively manufactured using the laser powder bed fusion (LPBF) process.

\begin{figure}
	\centering
	\includegraphics[keepaspectratio,width=0.5\textwidth]{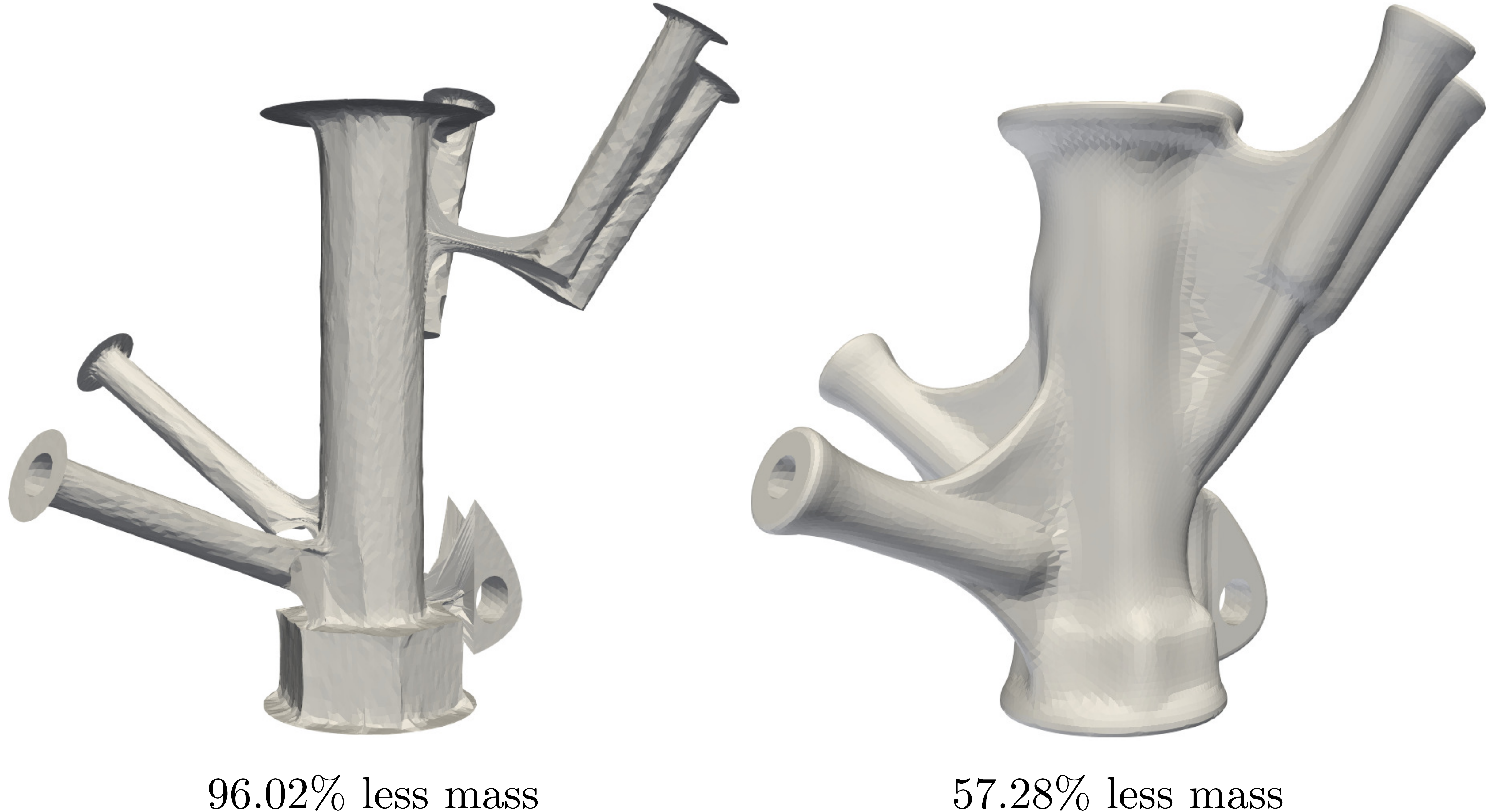}
	\caption{Unconstrained mass minimization of the considered tensegrity node using the bulk-surface Helmholtz shape filtering with $\beta=1$ (left) and the explicit filtering (right).}
\label{fig:tens-explicit-implicit-unconst-mass-min}
\end{figure}

Here, an aluminum node of the tower is reconsidered for shape optimization. Fig. \ref{fig:tensegrity-tower}b illustrates the problem setup including the boundary conditions for structural analysis and design optimization. Volume mesh consists of 212689 tetrahedrons and the outer surface of the node colored in green is subject to optimization. In the first step, mass of the node is optimized without any constraint. This is especially challenging due to large compression that mesh undergoes to reach the global optimal shape, which is basically projection of the design surface onto the non-design surfaces. Fig. \ref{fig:tens-explicit-implicit-unconst-mass-min} depicts optimal shapes obtained with the bulk-surface Helmholtz filtering using $\beta=1$ and the Vertex-Morphing technique. While the former handles boundary and internal nodes implicitly and simultaneously, the later treats them sequentially by first applying the explicit filtering on the surface and then using a pseudo-structural model with the Jacobian-based stiffening to deform the internal mesh. Furthermore, distorted or collapsed mesh is used as the optimization stopping criterion. Although both methods leverage the Jacobian-based stiffening technique to improve the deformed internal mesh quality, with the bulk-surface Helmholtz filtering much bigger surface deformation and objective improvement ($96.02\%$ vs. $57.28\%$) could be achieved.  

\begin{figure}
	\centering
	\begin{tabular}{@{}c@{}}
		\includegraphics[width=\linewidth,keepaspectratio]{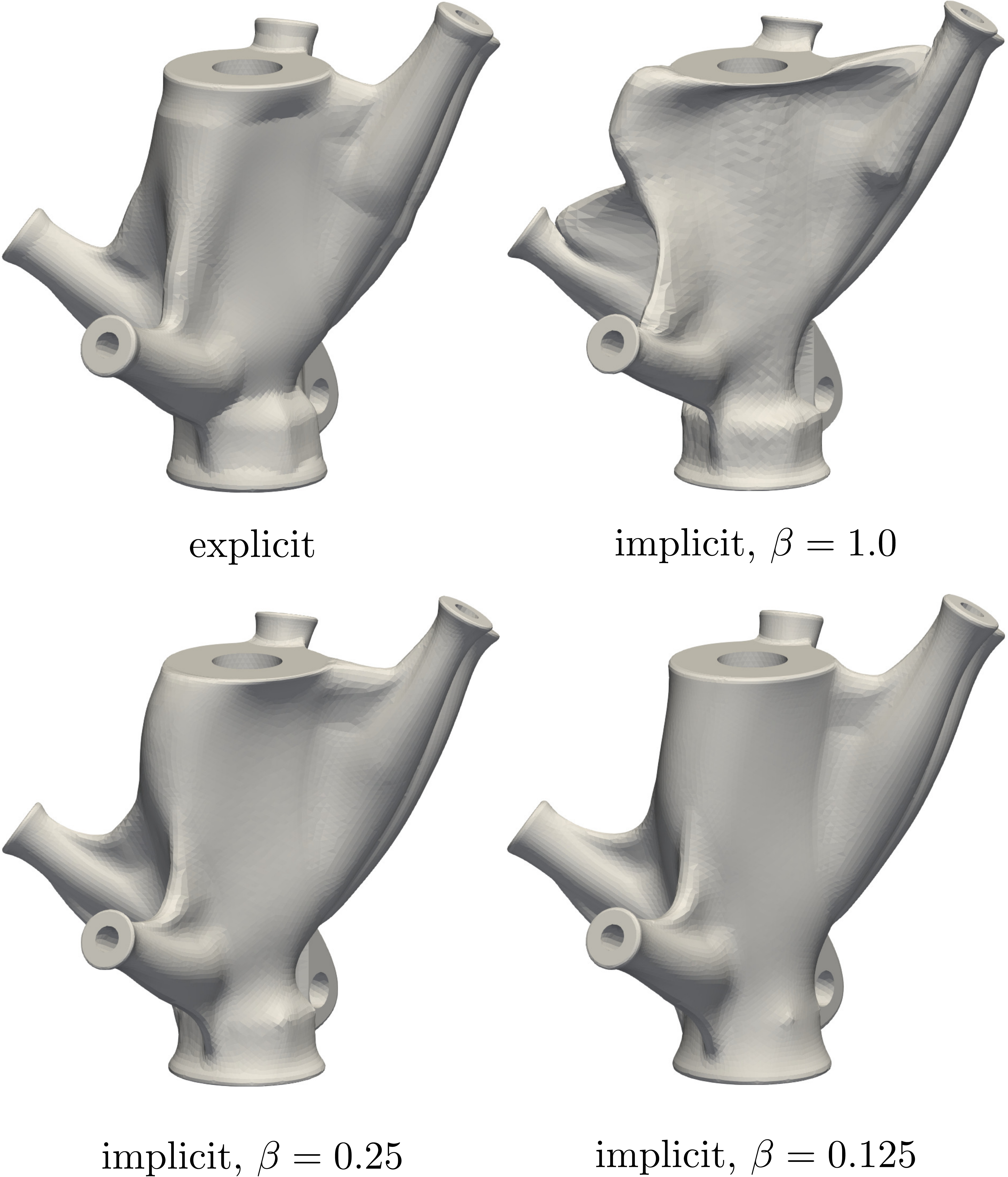} \\[\abovecaptionskip]
		\small (a) Geometries obtained with the explicit filter and the \\ implicit filter for various weighting factor $\beta$ .
	\end{tabular}
	
	\vspace{\floatsep}
	
	\begin{tabular}{@{}c@{}}
		\includegraphics[width=\linewidth,keepaspectratio]{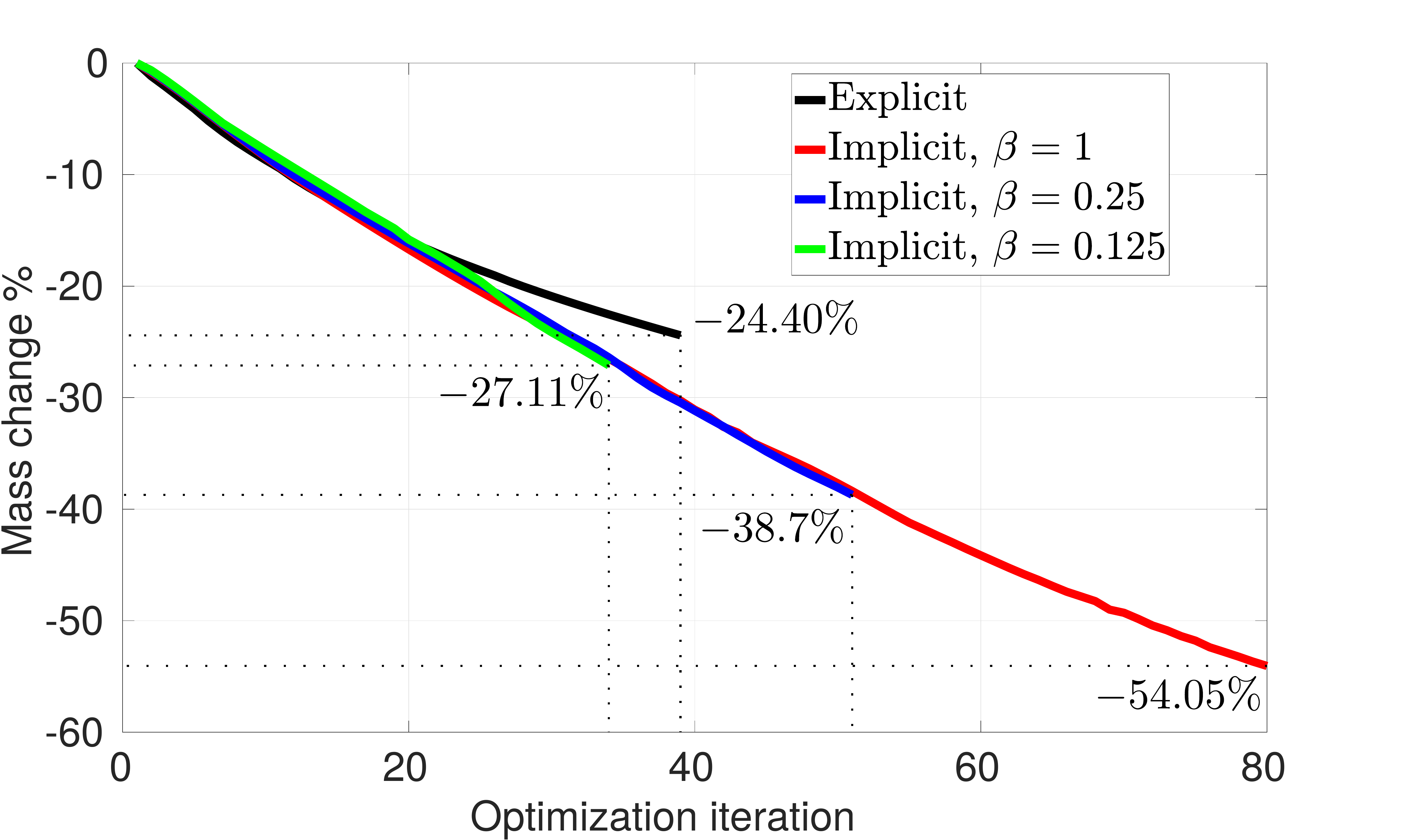} \\[\abovecaptionskip]
		\small (b) History of objective function until volume mesh distortion\\ (elements with negative Jacobian) stops FEM calculation.
	\end{tabular}
	
	\caption{Strain energy constrained mass minimization of the tensegrity node using the bulk-surface Helmholtz shape filtering (implicit) and the Vertex-Morphing technique (explicit).}
	\label{fig:tens-explicit-implicit-const-mass-min}
\end{figure}

\begin{figure}
	\centering
	\includegraphics[keepaspectratio,width=0.5\textwidth]{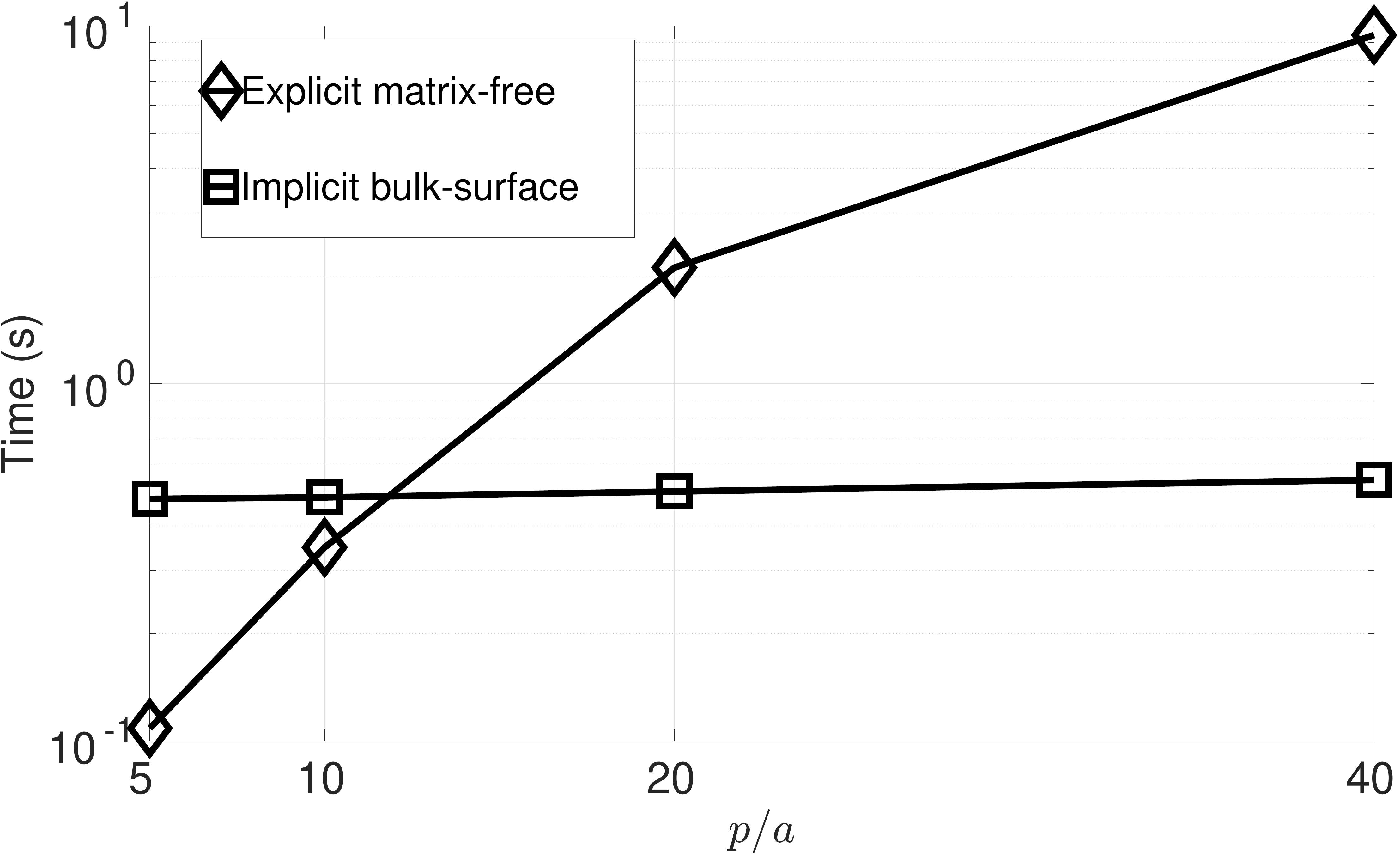}
	\caption{Time necessary for filtering a surface sensitivity field using the bulk-surface filter and the explicit surface filter without (matrix-free) storing the filter matrix for various $p/a$ ratios.}
	\label{fig:time-comp-knoten}
\end{figure}

To study the influence of the weighting factor $\beta$ on the performance of the proposed bulk-surface filtering, strain energy constrained mass minimization of the node is performed and the result are compared against that of the Vertex-Morphing technique. It should be noticed that the lower the factor is, the less attention or weight is given to the preservation of the internal mesh quality, which results in an earlier mesh distortion. This is clearly seen in the results presented in Fig. \ref{fig:tens-explicit-implicit-const-mass-min}. However, we observe that even for a small value $\beta=0.125$ the proposed implicit filtering outperforms the explicit one. On the other hand, we notice that the proposed bulk-surface filtering requires a linear system solve (Eq. \ref{e:implicit-sensitivity-mapping}b) for calculating each response function derivative with respect to solid control points. Comparison of the time necessary for filtering a surface sensitivity field using the explicit surface filter and the bulk-surface filter for different support span to element size ratios $p/a$ is shown in Fig. \ref{fig:time-comp-knoten}. Similar observations and conclusions as in the shell optimization case (Sec. \ref{sec:shell_opt}) can be made.  For moderate ratios, e.g. $p/a=10$, cost of the bulk-surface filter is comparable to that of the matrix-free explicit filtering, whereas for higher ratios cost of the matrix-free explicit filtering is hindering the performance.   

\section{Conclusions}
In this work, we visited two major filtering techniques for node-based shape optimization, namely the convolution-based (explicit) and the PDE-based (implicit). Consistency and mesh independency, as two strict requirements for any kind of filtering, were carefully discussed and demonstrated for both techniques. Mesh independent shape filtering was achieved by scaling discrete control sensitivities with the inverse of the mesh mass matrix. Supported by numerical experiments, a regularized Green's function was introduced as an equivalent explicit form of the so-called Helmholtz-/Sobolev-based (implicit) filter. Furthermore it was concluded that the Gaussian and linear kernel functions act as a low-pass filter suppressing fine scales, whereas the Helmholtz and regularized Green kernels are high-pass smoother, which allow generation of features like sharp edges and corners. With the aim of controlling the boundary smoothness and preserving the internal mesh quality simultaneously, this work introduced an implicit bulk-surface filtering technique for the shape optimization of volumetric domains. Its superior efficiency and robustness from the explicit filtering was demonstrated for shape optimization of a complex solid structure. Overall, whether surface or bulk domain is filtered, the implicit approach was found to be numerically more efficient and unconditionally consistent, compared to the explicit one.

\section*{Replication of results}
The software package used in this work, \emph{Kratos-Multiphysics}, is open-source and available for download at \url{https://github.com/KratosMultiphysics/Kratos}. The source code for replicating the results can be downloaded from \url{https://github.com/KratosMultiphysics/Kratos/tree/OptApp}. For the datasets that were used to generate the results, please contact the authors.

\section*{Conflict of interest}
On behalf of all authors, the corresponding author states that there is no conflict of interest. 
\section*{Acknowledgements}
This research is funded by Deutsche Forschungsgemeinschaft (DFG) through the project 414265976 TRR 277 C-02. The authors gratefully acknowledge the support.

\bibliographystyle{apalike}       

\bibliography{mybibfile}

\end{document}